\documentclass[preprint,12pt,authoryear]{article}
\usepackage{natbib}
\usepackage{amsmath,amssymb,amsfonts}
\usepackage{algorithmic}
\usepackage{graphicx}
\usepackage{textcomp}
\usepackage{hyperref}

\newtheorem{definition}{Definition}

\newtheorem{theorem}{Theorem}
\usepackage{booktabs}
\usepackage{tikz}
\usetikzlibrary{shapes,arrows}
\usetikzlibrary{positioning}
\tikzstyle{block} = [draw, fill=white, rectangle, 
    minimum height=3em, minimum width=3em]
\tikzstyle{circ} = [draw, fill=white, circle, minimum size=2.5em]
\tikzstyle{input} = [coordinate]
\tikzstyle{output} = [coordinate]
\tikzstyle{pinstyle} = [pin edge={to-,thin,black}]
\newcommand{\eq}{\overset{\mathrm{def}}{=}} 

\newcommand{\be}[0]{\begin{equation}}
\newcommand{\ee}[0]{\end{equation}}
\newcommand{\ben}[0]{\begin{equation*}}
\newcommand{\een}[0]{\end{equation*}}
\newcommand{\bena}[0]{\begin{eqnarray*}}
\newcommand{\eena}[0]{\end{eqnarray*}}
\newcommand{\bea}[0]{\begin{eqnarray}}
\newcommand{\eea}[0]{\end{eqnarray}}

\newcommand{\wh}{\widehat{\theta}}
\newcommand{\wt}{\widetilde{\theta}}
\newcommand{\RR}[0]{\mathbb R}
\newcommand{\tends}{\rightarrow}

\newcommand{\lm}[1]{\lim_{#1\tends\infty}}

\newcommand{\ltwo}{{\cal L}_2}

\newcommand{\bluec}[1]{\textcolor{black}{#1}}

\def\BibTeX{{\rm B\kern-.05em{\sc i\kern-.025em b}\kern-.08em
    T\kern-.1667em\lower.7ex\hbox{E}\kern-.125emX}}
\title{A Historical Perspective of Adaptive Control and Learning\thanks{The first author is supported by the Boeing Strategic University Initiative. The second author performed his part of the work in IPME RAS under support by Ministry of Science and Higher Education of the Russian Federation (Project no. 075- 15-2021-573 ). \\ Published at \url{https://doi.org/10.1016/j.arcontrol.2021.10.014}.}
}
\author{Anuradha M. Annaswamy and Alexander L. Fradkov  
}
\date{}
\begin{document}
\maketitle


\begin{abstract}
  This article provides a historical perspective of the field of adaptive control over the past seven decades and its intersection with learning. A chronology of key events over this large time-span, problem statements that the field has focused on, and key solutions are presented. Fundamental results related to \bluec{stability and robustness of adaptive systems and learning of unknown parameters} are sketched. A brief description of various applications of adaptive control reported over this period is included.
\end{abstract}
\section{Introduction}
The goal of adaptive control is real-time control of uncertain dynamic systems through adaptation and learning. This paper takes a historical perspective of the field of adaptive control over the past seven decades. Given the recent upsurge of interest in learning,  in the Machine Learning and Control communities,  both offline and online, such a perspective is timely and warranted. 

The scope that we aim to cover is clearly ambitious. Covering events that span 70 years, chronicled in more than 15 textbooks, 20 edited books, hundreds of surveys, and thousands of research publications in journals and conferences in 30 pages is a formidable task. The goal of this article is to accomplish this task by focusing on the highlights of this field,  emphasize key lessons learned,  delineate key solutions derived, and identify a few  takeaway messages. 

Here are the highlights of this article:
\begin{itemize}
\item Over the last seventy years, the field of adaptive control has witnessed advances in both deterministic and continuous-time systems and stochastic discrete-time systems. This article is one of the first attempts to trace the development in both  domains.
\item \bluec{The article focuses mainly on those advances in adaptive control that have a significant intersection with parameter learning.}
\item This article has made a concerted effort in chronicling key advances that have occurred globally.
\item	The article presents events chronologically (Section II), through problem statements (Section III), and through highlights of solutions (Section IV). An interested reader may delve into only one or more of these sections and do a deeper dive, if interested, by reading any of the 250 references listed at the end of the article.
\item	The article provides a snapshot of various applications of adaptive control (Section V).
\item	Concluding remarks including a few takeaway messages  are provided in Section VI.
\end{itemize}

When it comes to real-time control of uncertain dynamic systems, the efforts of the control community extend significantly beyond adaptive control. There are several topics that are at the boundaries, such as sliding-mode control, iterative learning control, and linear-parameter-varying control that are not addressed in this survey. While we devote a section to the topic of reinforcement learning (RL) that has a rich intersection and complementarity to adaptive control, we defer the reader to other articles  for a more  comprehensive discussion on RL as well as other topics that lie at the intersection of adaptation and learning.

\section{A chronology}
The history of adaptive control systems is almost as long as the entire field of control systems, as the concept of adaptation is fairly close to the notion of feedback. As such, this concept has been explored from the 1950s to the present and continues to be an area of intense activity. We classify various developments in this area into three chunks of time, 1950-65, 1965-85, and 1990s-\bluec{2000}. During the first 15 years, several contributions arose separately in the context of deterministic continuous-time systems and pattern-recognition (which employed a stochastic framework), and are organized under these two broad headings. Subsequent sections outline parallel developments in deterministic and stochastic systems in a combined manner. While the 70s and 80s witnessed the development of a stability framework, later decades developed a robustness framework for the adaptive systems. Key developments in all of these decades are outlined below.

\subsection{1950-65}
\subsubsection{Deterministic and continuous time}
The term adaptation is defined in biology as “an advantageous conformation of an organism to changes in its environment.” The earliest reflection of this principle in an engineering context can be found in \citep{Dre57}\footnote{Origins of adaptation rules  can be traced even earlier to 1949, in the form of Hebbian rules \citep{Hebb49} that connected weight adjustments in a neuron to performance.}. The authors coopted this fundamental principle in their definition of an adaptive system in the context of a control system, and defined an adaptive control system to be one which monitored its own performance and adjusted its parameters in the direction of better performance \cite{Dre57}. The implicit implication here is that a non-adaptive system would then have parameters that are fixed and not adjusted. To provide more clarity, and distinguish an adaptive system from a non-adaptive one, references \cite{Aseltine58} and \cite{Stromer59} introduced definitions of adaptive systems. In fact, there was a profusion of definitions of adaptive systems at this time based on what was adapted, what  the adaptation was in response to, time-scales of adaptation, or from whose viewpoint. It could be argued that the classes of adaptive systems outlined in \citep{Aseltine58} are precursors to the current approaches in adaptive control. 

Similar to \cite{Dre57}, the authors of \cite{Whit58} focused on a servo problem where the process output was required to follow a commanded output in the presence of parametric uncertainties.  They developed what came to be known as the {\it MIT-rule} as a core adaptive mechanism which served as an outer loop with the inner loop consisting of a standard feedback control system. The adaptive mechanism adjusted the control parameter $\theta$ using a simple rule
\be \dot \theta = - k e(t) \nabla_\theta e(t) \label{mit-rule}\ee
where $e(t)$ denoted a tracking error between the process output $y(t)$ and a reference output $y_m(t)$, and $\nabla$ stands for the gradient. The idea therefore is to have the adaptive mechanism use \eqref{mit-rule} to \bluec{estimate}
, i.e. learn the correct value of the parameter that the feedback controller in the inner loop must deploy. \bluec{This main idea continues to pervade all adaptive control methods to-date.}

The motivation for the study and implementation of adaptive control systems came from applications in aerospace – for autopilot design in flight control  (\cite{Greg59,Ham13}). As high performance aircraft  routinely encounter a wide range of operating conditions, there was a need to develop sophisticated regulators that would adapt their parameters online so that they are not constrained to work with constant gains which may limit their operation to a small flight envelope. This led to several symposia on adaptive systems in the early 60s, with what was referred to as a {\it three-legged milking stool} for advanced flight control systems that consisted of aerodynamics, GNC (Guidance, Navigation, and Control), and adaptation \citep{Ham13}. Around the same time, Bellman and Kalaba introduced the term adaptive in the context of multistage decision processes as belonging to the last of a series of three stages in the evolution of control processes. With the first two denoted as deterministic and stochastic control processes, an adaptive control process was defined as when the designer has very little knowledge about the system dynamics or even the statistics of any random inputs that may be present (\cite{Bel59,Bel61}). Yet another early evidence of interest in adaptive regulators is a patent by Caldwell \citep{Cald50}. Several additional references can be found in \cite{Astrom96} and in Chapter 1 in \cite{Narendra2005}.

Cautionary inputs and guidance for the design of the controllers and the adaptive mechanisms soon started to appear. Any successful adaptive system has to cope with changes in its environment for its survival and performance. Familiarity with the environment results in better understanding, and better understanding enables the system to better predict the changes in the system. However, understanding and controlling are two distinct activities; predictive ability does not translate directly into the ability to control. Often the converse may be true – a better ability to control may help in a better understanding but here too there may be limits. Needless to say, the connection between identification and control is complex and was explored in a number of seminal papers and textbooks during the '60s. One of them is Feldbaum’s concept of dual control \citep{Feldbaum60}, that emphasized the need for an optimal control action that is taken for a system with uncertainties. Feldbaum pointed out that the requisite control has to have dual components, one of probing for enhancing identification and one of caution for ensuring stable control action. Too much of a focus on identification may not result in satisfactory control; too much emphasis on controlling the system may not lead to \bluec{satisfactory} learning. The design of dual control with the right mix of both of these components is therefore a huge challenge and the grand goal of the field of adaptive control. These  two intertwined concepts of identification and control pervade Machine Learning (ML) as well (\cite{Kaelbling96,ishii2002}), and often go under the monikers of “exploration” and “exploitation.” 


\subsubsection{Pattern Recognition and Classification}
A parallel development of adaptation can be traced in the field of pattern recognition and classification, which occurred during the same period. As the title of \cite{Widrow64} attests, it was observed that a gradient descent type algorithm, similar to that in \eqref{mit-rule}, plays a central role not only in control problems but also in pattern recognition. In addition to \citep{Widrow61,Widrow64,Abramson63}, several groups in USSR led by Aizerman \cite{Aizerman63}, Lerner \cite{Vapnik63,Vapnik64}, Yakubovich \cite{Yakubovich63,Yakubovich65}, and others \citep{Bongard61,Braverman62,Fradkov20,FradkovPolyak20} developed deterministic and stochastic approaches for input classification and pattern recognition.  A common element to a diverse set of problems in pattern recognition, signal processing \citep{Widrow67}, and perceptrons \citep{Rosenblatt62}, was the determination of a set of parameters or weights that leads to desired classification, filtering, or tracking performance using input-output data. In contrast to the earlier discussions of control systems, the treatment in these works was in discrete-time rather than continuous time, and instead of a deterministic framework, employed a stochastic framework with noisy measurements and inputs. Widrow’s Adaline filters also led to the foundation of neural networks, deep and otherwise \citep{Widrow1990}. The approach taken in many of these works was statistical in nature, with their foundations in communications and decision theory \citep{Marill60,Widrow59}. Except for brief mentions, this survey will not focus on the evolution of pattern recognition or its intersection with adaptive control.
\subsection{1965-1985}

It was soon realized that the MIT-rule proposed in \cite{Whit58} can result in instability, especially when there is sufficient phase lag between the measurement of error and adjustment of the parameters. Several authors contributed to the formulation of a stability framework for the analysis and synthesis of adaptive systems where real-time decisions in the form of parameter adjustment in dynamic systems were taken using online data. Notable ones came from the authors of \cite{Gray63}, \cite{Shack65}, \cite{Parks66}, \cite{Mon67}, and \cite{Nar74}.   Lyapunov’s method was suggested in lieu of a gradient descent approach as in \eqref{mit-rule}, and ended up as the foundation for stability of adaptive systems\footnote{In hindsight, the MIT-rule can be viewed as a partial Lyapunov function, as it only included an $\ltwo$-norm of the performance error in its cost function.}. Independently, the same problem with similar conceptual tradeoffs was also addressed in  deterministic discrete time setting  by Yakubovich in \cite{Yak68,Yakubovich72}. Several seminal results were published during this period which witnessed  surveys by Lindorff and Carroll \citep{Lind73}, Landau \citep{Landau74}, Wittenmark \citep{Wit75}, Unbehauen \citep{Unb75}, and others \citep{Ash76,Parks80,Voronov84}. These were followed by edited books such as \cite{Nar80,Unb80,Har80}, and subsequent textbooks in deterministic and continuous-time \cite{Fom81,Narendra_1989,Astrom_1995,Ioannou1996,Sastry_1989,Tao03,Krstic_1995,Fradkov99}, stochastic systems in books and papers such as \cite{Kumar1986,Duncan90,Borkar79,Becker85,Fomin91b}, and multiple-input, multiple-output systems in \cite{Tao03}. These addressed adaptive control architectures and algorithms for a range of dynamic systems, with either full-state or partial-state measurements available in real-time. The efforts during these 15 years laid the foundation for stable adaptation in dynamic systems, both deterministic and stochastic, where the uncertainties were predominantly in their parameters. The overall goal was to ensure a closed-loop system that was well-behaved and met control goals such as tracking and regulation asymptotically.

In deterministic systems, the  structure of the algorithm for adjusting their parameter $\theta$ was of the form
\be \dot \theta(t) = - k(t) e(t) \phi(t)  \label{stability-based}\ee
where $\phi$ is a suitably chosen regressor that may or may not coincide with the gradient of a well-defined loss function, and $k(t)$ represents a normalization component. The choices of $k$ and $\phi$ were guided by the determination of an underlying Lyapunov function and the interplay between the adjustable parameters and the signals in the closed-loop system, leading to an approach that is most commonly  termed \textit{Model Reference Adaptive Control (MRAC)} and used in deterministic continuous-time systems.
In stochastic systems, the works by Astrom and coworkers \citep{Astrom73,Astrom_1995} led to Self-tuning Regulators (STR) associated with minimum variance controller with their foundation laid in papers such as \citep{Lju77b},\citep{Solo79,Landau82,Bit83,Kumar83,clarke1985generalized,Johan96}.
In all these cases, conditions under which learning, that is, accurate parameter estimation, can take place were precisely articulated. Both necessary and sufficient conditions were derived \citep{Morgan_1977,Lju77a,Lju77b,Lju83,And85}.

Yet another link between adaptation and learning is due to Yakov Tsypkin who proposed a unified framework based on stochastic approximation machinery. Parameter choice and convergence results then follow from the results on stochastic approximation obtained earlier in mathematical statistics based on average risk  minimization \citep{Tsypkin66,Tsypkin71}. We defer the details of the problem statement to Section III. 

As evidenced by the chronology above, there are two parallel, and very similar evolution of the branches of adaptive control in deterministic systems and stochastic systems, with the obvious distinction associated with the underlying tools. The problem of convergence of the tracking error in the former case had a counter-part of a minimum variance controller in the latter. The term adaptive controller remained in vogue for deterministic systems and its counterpart in stochastic systems was termed self-tuning regulators; the terms adaptation and self-tuning were used synonymously. The fundamental tenets of stability and convergence in adaptive systems and tradeoffs between performance and learning were however found to be invariant to these two branches. \bluec{We note that for stochastic systems, our focus in this paper is restricted to adaptation and parameter learning in  discrete-time systems. There is a significant and rich literature present in adaptive control of stochastic continuous-time  systems as well (see \cite{Wertz89,Gever91,Caines92,Duncan99} for linear systems  and \cite{LiKrstic20} for nonlinear systems). } \bluec{Most of these ideas and results have discrete-time counterparts, which are presented in brief in the following sections.} Details of the problem statements are postponed until section III.

\subsection{1990s-present}
\subsubsection{Adaptive Control of Deterministic and Stochastic Systems}
With the stability framework established in the 70s, the next broad milestone in the evolution of adaptive control systems was a robustness framework established in the 80s with textbooks capturing the details of various solutions in the 90s. It was soon realized that both gradient algorithms as in \eqref{mit-rule} and stability-based algorithms that employed a Lyapunov approach as in \eqref{stability-based} were inadequate in ensuring robustness to perturbations such as bounded disturbances and unmodeled dynamics \citep{Roh82}. Several approaches were developed around the same time \citep{Ega79,Pet82,Kre82,Ioa84,Praly83,Praly84,And86,Nar86,Nar87,Ioa86,Mid88,Tsa89,Ort89,Naik92} in ensuring that adaptive control systems not only provided appropriate adaptation to accommodate parametric uncertainties but also provided robustness to withstand non-parametric uncertainties such as external disturbances, time-varying parameters, and unmodeled dynamics.  Broadly, these approaches either relied on properties of persistent excitation of the exogenous command signals \citep{And86,Nar86} with the same adaptive laws as in \eqref{stability-based}, or in modifying the adaptive law in a suitable manner \citep{Fradkov79,Pet82,Kre82,Ioa86,Nar87,Fradkov87,Wen92}. A parallel to the latter corresponds to the use of regularization in machine learning \citep{Gaudio20AC}. Similar results can be found in discrete time as well (for example, \citep{disIoa86,Tao95,Wen92,Clu88}). \bluec{The use of diminished persistent excitation with time was utilized to obtain an elegant framework for adaptive optimal control in stochastic systems in \citep{GuoChen91,Guo95,Duncan99}.}

The above stability and robustness arguments also set the stage for addressing the control of nonlinear systems with parametric uncertainties. This too was addressed starting in the 90s, spawning a huge area of research with dozens of researchers laying the foundation of key results (see for example \cite{Krstic_1995}). Special classes of adaptive nonlinear systems that arise in robotics were addressed at length in \cite{Slotine1991} even earlier. Methods based on feedback linearization, backstepping, and averaging led to several seminal results in this area. 
A class of problems related to control of nonlinear systems using neural networks witnessed significant activity during this period as well (see for example, \cite{Narendra_1991,Sanner_1992,Rov94,Polycarpou96,Yu_1996,Yu_1998,Hov08,Ren10}).

\subsubsection{Reinforcement learning/Approximate Dynamic Programming}
Towards the end of the 1980s the approach of reinforcement learning \citep{Sut18,Kaelbling96} 
was formulated and in the early 1990s  strong ties were identified between these topics and adaptive optimal control. A case in point is the reference \cite{Sutton_1992}, clearly indicated in its title, \textit{Reinforcement learning is direct adaptive optimal control.}
In the works on control based on RL, a performance index is introduced, usually as an integral functional, and neural networks are used to approximate either  the predicted optimal value (Bellman function) of this functional (known as Value Iteration (VI)) or the optimal control policy (policy iteration), based on the HJB equation. Such an approach is referred to as  ``approximate dynamic programming", "neuro dynamic programming" \citep{Pow07,Bert08}, or ``adaptive dynamic programming" \citep{Lew09}. 
Analytical frameworks for cases when the state and action sets are finite and for the more difficult case when they are infinite have been addressed. Related problem statements are briefly addressed in Section III.


\section{Problem Statements}
This section outlines the problem statements that have been proposed under the rubric of adaptive control. We classify them into four categories, the first three of which are based on whether they are in continuous-time or discrete-time, and deterministic or stochastic. We do not address stochastic continuous-time systems in this paper mostly since that area was developed much in parallel with the discrete-time case \citep{Wertz89,Gever91,Caines92,Duncan99,LiKrstic20}. A brief discussion on RL is also included in this section. 

\subsection{Deterministic and Continuous-time Systems}\label{determ-problem}
The aim in adaptive control problems is to design an exogenous input $u(t)\in\RR^m$ that affects the dynamics of a system given by \bea \dot x &=& f(x,\theta,u,t)\nonumber\\
y&=& g(x,\theta,u,t)\label{plant1}\eea
where $x(t)\in \RR^n$ represents the system state, $y(t)\in\RR^p$ represents all measurable system outputs, with many physical systems obeying the inequality $n>>p>m$ \citep{Qu2020}. $\theta\in\RR^\ell$ represents system parameters that may be unknown, and $f(\cdot)$ and $g(\cdot)$ denote system dynamics, that may be nonlinear, that capture the underlying physics of the system.  The  functions $f(\cdot)$ and $g(\cdot)$ also vary with $t$, as disturbances (often modeled as deterministic quantities) and stochastic noise may affect the states and output. The goal is to choose $u(t)$ so that $y(t)$ tracks a desired command signal $y_c(t)$ at all $t$, and so that an underlying cost $J((y-y_c),x,u)$ is minimized.  
In what follows, we will refer to the system that is being controlled as a plant. 

As the description of the system as in \eqref{plant1} is based on a plant model, and as the goal is to determine the control input in real time, all control approaches  make  assumptions regarding what is known and unknown. To begin with, as the plant is subject to various perturbations and modeling errors due to environmental changes, complexities in the underlying mechanisms, aging, and anomalies, both $f$ and $g$ 
are not fully known. The field of adaptive control has taken a parametric approach to distinguish the known parts from the unknown. In particular, it is assumed that $f$ is a known function, while the parameter $\theta$ is unknown. A real time control input is then designed so as to ensure that the  tracking goals are achieved by including an adaptive component that attempts to estimate the parameters online.  A linearized version of the problem in \eqref{plant1} is of the form 
\be y=W(s,\theta)[u] \label{plant-linear}\ee
where $s$ denotes the differential operator $d/dt$, $W(s,\cdot)$ is a rational operator of $s$, and $\theta$ is an unknown parameter, and the goals of tracking and regulation are the same as above.

In the following  subsections, four broad categories of subproblems that have been addressed in the context of adaptive control in deterministic continuous-time systems are described.

\subsubsection{Boundedness and real-time decision making}
As mentioned above, the control goal is to ensure that 
\be \lm{t} e(t)=0 \label{adapt-goal} \ee where $e(t)=y(t)-y_c(t)$. As these decisions are required to be made in real time, the focus of the solutions  is to have them lead to a closed-loop dynamic system that has bounded solutions at all time $t$ and a desired asymptotic behavior. The central question, therefore, is if this can be ensured even when there are parametric uncertainties in $\theta$ and several other non-parametric uncertainties that may due to unmodeled dynamics, disturbances, and the like. Once this is guaranteed, the question of learning, in the form of parameter convergence, is addressed. As a result, {\it control for learning} is a central question that is pursued in the class of problems addressed in adaptive control rather than \textit{learning for control} \citep{Krs21}.

\subsubsection{Certainty Equivalence Principle and Adaptive Control Solutions}\label{ss:cep}
The well known certainty equivalence principle (CEP) proceeds with the following mantra: first, optimize under perfect foresight, then substitute optimal estimates for unknown  values. This philosophy underlies  all adaptive control solutions by first determining a controller structure that leads to an optimal solution when the parameters are known and then replace the parameters in the controller with their estimates. The difficulty in adopting this philosophy to its fullest stems from the dual nature of the adaptive controller, as it attempts to accomplish two tasks simultaneously, estimation and control. This simultaneous action introduces a strong nonlinearity into the picture and therefore renders a true deployment of the certainty equivalence principle difficult if not impossible. The procedure for adaptive control is therefore modified, with the first step corresponding to a controller that leads to a \textit{stable} solution rather than optimal one. In other words, much of the adaptive control literature has focused on deriving stable solutions first and foremost for the real-time control of systems with parametric uncertainties, followed by an effort to estimate the unknown parameters, and optimization addressed at the final step. Such a breakdown of the problem overcomes the intractability of the certainty equivalence principle and leads to tractable procedures.

A typical solution of the adaptive controller takes the form
\bea u&=& C_1(\theta_c(t),\phi(t),t) \label{ad1} \\\dot\theta_c &=& C_2(\theta_c,\phi,t) \label{ad2} \eea
where $\theta_c(t)$ is an estimate of a control parameter that is intentionally varied as a function of time, $\phi(t)$ represents all available data at time $t$.  The nonautonomous nature of $C_1$ $C_2$ is due to the presence of exogenous signals such as set points and command signals. The functions $C_1(\cdot)$ and $C_2(\cdot)$ are deterministic constructions, and make the overall closed-loop system nonlinear and nonautonomous. 
The challenge is to suitably construct functions $C_1(t)$ and $C_2(t)$ so as to have $\theta_c(t)$ learn the requisite unknown control parameter $\theta_c^*$, and ensure that stability and asymptotic stability properties of the overall adaptive systems are ensured. These constructions have been delineated for deterministic systems in \citep{Narendra_1989,Astrom_1995,Ioannou1996,Sastry_1989,Krstic_1995,Tao03} and other textbooks. The solutions in these books and several papers in premier control journals such as Transactions on Automatic Control and Automatica have laid the foundation for the construction of $C_1$ and $C_2$ for a large class of dynamic systems in \eqref{plant1}. 

\noindent\underline{Model Reference Adaptive Control}
A tractable procedure for determining the structure of the functions $C_1$ and $C_2$, denoted as \textit{Model Reference Adaptive Control}, uses the notion of a reference model, and a two-step design consisting of an algebraic part for determining $C_1$ and an analytic part for finding $C_2$. A reference model provides a structure to the class of command signals $y_c(t)$ that the plant output $y$ can follow. For a controller to exist for a given plant-model using which the closed-loop system can guarantee output following, the signal $y_c$ needs to be constrained in some sense. A reference model is introduced to provide such a constraint. In particular, a model ${\cal M}$ and a reference input $r$ is designed in such a way that the output $y_m(t)$ of ${\cal M}$ for an input $r(t)$ approximates the class of signals $y_c(t)$ that is desired to be followed. With a reference model in ${\cal M}$, the algebraic part of the MRAC corresponds to the choice of $C_1$ with a fixed parameter $\theta_c^*$ such that if $\theta_c(t)\equiv\theta_c^*$ in \eqref{ad1}, then $\lm{t}y_p(t)-y_m(t)=0$. With such a $C_1$ determined, noting that $\theta_c^*$ could be unknown due to the parameteric uncertainty in the plant, the analytic part focuses on finding $C_2$ such that output following takes place with the closed-loop system remaining bounded.

\subsubsection{Learning $=$ parameter estimation}
With the problem statement as above, it is perhaps clear to the reader that the organic connection between the adaptive control problem and learning enters through parameters. Given that what's unknown about the dynamics is the plant parameter $\theta$, or equivalently the control parameter $\theta^*_c$, learning is synonymous with accurate parameter estimation. That is, it is of interest to have the parameter estimate $\theta_c$ converge to $\theta_c^*$ in the context of a control problem, and in identification problems, for an estimate $\wh$ to converge to $\theta$. The goal in either case is to determine conditions under which this convergence take place. These conditions are linked to  properties defined as persistent excitation (PE) and uniform observability \citep{Narendra_1987,Boyd_1983,Morgan_1977,Anderson_1982,Jenkins_2019}. These PE properties are usually associated with the underlying regressor $\phi$, and typically realized by choosing the exogenous signals such as $r(t)$, the input into the reference model ${\cal M}$ appropriately, which the control designers have the freedom to select. 

\bluec{Yet another extension that has been explored successfully in adaptive control is the notion of multiple models \citep{narendra1997multiple}. The goal is  the same as in MRAC, but to  accomplish adaptation rapidly. As the name suggests, the solution consists of generating multiple models of the plant, with multiple identification errors, one associated with each model, and carry out two steps, of switching and tuning. Switching consists of determining the model with the smallest error using a suitable criterion, and tuning corresponds to the adjustment rule that identifies the parameters of that particular model. Several algorithms are suggested in \citep{narendra1997multiple} and the references therein. In \cite{narendra1997multiple}, the premise is that $p^*$, the plant parameter suddenly changes, and the goal is to quickly determine an adaptive controller using a combination of fixed models and adaptive models where the plant parameter is identified. While learning is a part of the objective of adaptive models, the focus of the paper is primarily in determining a closed-loop system that remains stable. The counterpart to the concept of multiple model-based adaptive control in the fixed control domain is supervisory control \citep{morse1996supervisory}.}

\subsubsection{Robust adaptive control}
The assumption that the uncertainties in \eqref{plant1} and \eqref{plant-linear} are limited to just the parameter $\theta$, and that otherwise $f$ and $g$ or $W(s)$ are known, is indeed an idealization. Several departures from this assumption can take place in the form of unmodeled dynamics, time-varying parameters, disturbances, and noise. For example, the linear plant may have a form
\be y=\left[W(s,(\theta(t)))+\Delta(s)\right][u+d(t)+n(t)] \label{plant-linear-robust}\ee where $d(t)$ is an exogenous bounded disturbance, $n(t)$ represents measurement noise, the parameter $\theta$ is time-varying and is of the form
\be \theta(t)=\theta^*+\vartheta(t)\label{timevarying-parameter}\ee where $\theta^*$ is an unknown constant parameter but is accompanied by additional unknown variations in the form of $\vartheta(t)$, and $\Delta(s)$ may represent higher-order dynamics that is either not known, poorly known, or even deliberately ignored for the sake of computational simplicity. In all of these cases, a robust adaptive controller needs to be designed  to ensure that the underlying signals remain bounded, with errors that are proportional to the size of these perturbations. As mentioned earlier, 
these approaches either relied on properties of persistent excitation of the exogenous command signals \citep{And86,Nar86} with the same adaptive laws as in \eqref{stability-based}, or in modifying the adaptive law in a suitable manner \citep{Fradkov79,Pet82,Kre82,Ioa86,Nar87,Fradkov87,Wen92}. These are summarized in \citep{Narendra2005,Astrom_1995,Ioannou1996,Sastry_1989,Tao03,Krstic_1995,Fradkov99}. Details of these approaches are deferred to the next section.

\subsection{Stochastic and Discrete-time Systems}
A parallel development in adaptive control is one where the control decisions take place in a stochastic environment. The problem statements once again center around systems that are not known, with a random or noisy behavior being an essential feature. Here too, there are multiple classes of problems that have been studied over the past five decades, a broad division corresponding to Bayesian and Non-Bayesian problem statements \citep{Kumar1986}. In both classes, similar to the problem statement in Section \ref{determ-problem}, the unknown part of the system pertains to its parameters. The former corresponds, as the name suggests, to problems where a  probability distribution of the parameter is known {\it a priori}, while in the latter, only a known set $\Theta$ to which the parameter belongs is given. Examples of the former include the Bayesian N-armed bandit problem \citep{Git79,Kumar81,Kumar85}, and self-tuning regulators (STR) \citep{Astrom_1995} for the latter. In this paper, we limit our discussion to the latter. 
\subsubsection{Self-tuning Regulators}\label{sec:STR}
The starting point for the STR problem is a Nonlinear Auto-Regressive Moving-Average model with noise (NARMAX) of the form
\begin{equation}
  \label{eq:10}
  \begin{split}
    y_k = &\sum_{i=1}^{n}a_{i}^*y_{k-i} + \sum_{j=1}^{m}b_{j}^*u_{k-j-d} +\sum_{i=0}^{n}c_{i}^*w_{k-i}\\
    &+ \sum_{\ell=1}^pd^*_{\ell}f_\ell(y_{k-1}, \ldots, y_{k-n}, u_{k-1-d}, \ldots, u_{k-m-d}),
  \end{split}
\end{equation}
where $a_{i}^*$, $b_{i}^*$,  $c_{i}^*$, and $d_i^*$ are unknown parameters  and $d$ is a known time-delay. $w_i$ corresponds to a white noise, stemming from measurement noise as well as input noise. The function $f_\ell$ is an analytic function of its arguments and is assumed to be such that the system in \eqref{eq:10} is bounded-input-bounded-output (BIBO) stable. The NARMAX plant in \eqref{eq:10} reduces to a linear ARMAX-model when the last term on the right-hand-side is absent. Alternate state-space representations rather than the input-output one in \eqref{eq:10} have also been analyzed in the literature. The goal once again is to determine the control input $u_k$, in real-time, so that the output $y_t$ either tracks a desired signal $y_{d,t}$ or is regulated around zero.

Similar to the deterministic counterpart, the adaptive control structures that solve the problem for stochastic systems have also centered around one that leads to a solution when the parameters are known.  We rewrite the system in \eqref{eq:10}, assuming that there are no nonlinearities, as 
\be A(z) y_{k+1}=z^d B(z) u_{k+1}+C(z)w_{k+1}\label{kumar1}\ee
where $A(z)=1-a_0z-\ldots +a_nz^{n+1}$,
$B(z)=b_0+b_1z +\ldots+b_pz^p$, 
$C(z)=1+c_0z+\ldots c_nz^{n+1}$, and $z$ denotes the shift operator, $zy_k=y_{k-1}$. Defining polynomials $F(z)$ and $G(z)$ as
$F(z)=1+f_1z+\ldots+f_{d-1}z^{d-1}$,
$G(z)=-g_0+g_1z+\ldots +g_nz^{n}$, we can express the polynomial $C$ as
\be C=AF+z^dG\label{kumar2}\ee
where we have suppressed the arguments for ease of exposition. $F(\cdot)$ represents the quotient of $A$ with respect to $C$. With the introduction of $F$ and $G$, we can rewrite the system in \eqref{kumar1} as 
\be y_{k+1}=z^d\left(\frac{B}{A}u_{k+1}+\frac{G}{A}w_{k+1}\right)+Fw_{k+1}\label{kumar3}\ee
It is easy to see that the desired control input is given by
\be u_{k}=-\frac{G}{BF}y_k\eq C(\phi_k,\theta_c^*)\label{kumar4}\ee
where $\phi_k = [y_{k-1}, \ldots, y_{k-n}, u_{k-1-d}, \ldots, u_{k-m-d}]^\top$. The parameter $\theta_c^*$ is a transformation of the parameters of $A$ and $B$, by virtue of the relation in \eqref{kumar2}. It can also be shown that \citep{Kumar1986} that the control input in \eqref{kumar4} minimizes the variance $E((1/N)\sum_1^N(y_k^2))$, and is often referred to as a minimum variance control  \citep{Astrom_1995,clarke1985generalized,Johan96}.

The self-tuning regulator addresses the design of a minimum variance control when the parameter $\theta^*$ is unknown. The corresponding solution pertains to the choice of the control input of the form  \citep{Ast70} 
 \be u_k=C(\phi_k,\widehat\theta_k)\label{kumar5}\ee
 and finding parameter updates for the parameter estimate $\widehat \theta$ so that minimum variance can be  achieved. We defer a discussion of various results related to STR to Section IV.
 
\subsubsection{Parameter Estimation and Persistent Excitation}\label{sec:PE}
Similar to the deterministic case discussed above, here too learning is tied with estimation of unknown parameters. The ARMAX problem considered in \eqref{kumar1} can be rewritten as a linear regression equation
\begin{equation}
  \label{eq:11}
  y_k = \phi_{k-1}^\top\theta_k^*+v_k,
\end{equation}
where $v_k$ is a noise term, not necessarily white, and $\theta^*$ is a vector of unknown parameters that needs to be estimated. Parameter estimation can then be carried out using a variety of iterative algorithms such as stochastic approximation \citep{Kumar83,Goodwin_1981} also known as stochastic gradient descent (SGD), and  recursive least squares \citep{Kumar85,Goodwin_1981}. As will be seen in Section IV, the conditions under which the estimates generated by these algorithms converge to the true values are well understood, also denoted as persistent excitation. The same procedure can also be adopted in adaptive control by starting with \eqref{kumar5}, and noting that it can be expressed once again as a linear regression. We discuss these details in Section IV.

A more generic formulation of the adaptation and learning problem was proposed by Tsypkin \citep{Tsypkin66,Tsypkin71} based on minimization of a averaged performance index. That is, the problem was posed as
$\min_{\theta} J(\theta)$  where $J(\theta)$ is the average of the cost function  $Q(x,\theta)$ over $x$ with an unknown density $p(x)$, where $x$ is the state and $\theta$ is a decision variable:
\be J(\theta)=\int_XQ(x,\theta)p(x)dx=E_xQ(x,\theta).\label{e:bellman}\ee
Tsypkin proposed a solution based on SGD as
$$\theta[n]=\theta[n-1]-\gamma[n]\nabla Q(x[n],\theta[n-1])$$
 where $\theta[n]$ is a recursive estimate of $\theta$. Choosing the cost function in an appropriate way allowed the author to present  different classes of algorithms described previously in the literature and a number of new ones  in a unified manner. He also showed that convergence of the algorithms can be established based on the stochastic approximation scheme under conditions of convexity and bounded growth of $J(\theta)$ and classical Robbins-Monro conditions on $\gamma[n]$ \citep{Rob52}, namely
 \be
\gamma[n]>0, \sum_n\gamma[n]=\infty, \sum_n\gamma^2[n]<\infty.\label{robbinsmunro}\ee

\subsubsection{Adaptive Optimal Control of Linear Quadratic Gaussian Systems}
The problem statements in sections \ref{determ-problem}, \ref{sec:STR}, and \ref{sec:PE} have focused on ensuring that a tracking error in states, or an output variance in the context of regulation is minimized \citep{Astrom_1995,Cla79}. An alternate class of problems has focused on minimizing a quadratic cost not only in states but also in the inputs. A typical problem formulation in this class is of the form \citep{Becker85}
\be x_{k+1}=Ax_k+B u_k+w_{k+1}\label{kumar10}\ee
where $A$ and $B$ are unknown matrices, and $w_k$ is a noise process made up of Gaussian i.i.d. random variables $N(0,1)$. The control objective is to determine $u_k$ such that the cost function 
\be J(A,B)\eq \limsup_{T\rightarrow\infty}\frac{1}{T}\sum_{i=1}^T\left[x_i^TQx_i+u_i^TRu_i\right], \;\; \label{kumar11}\ee
where $Q=Q^T>0, R=R^T>0$, is minimized. 


\subsection{Deterministic and Discrete-time Systems}\label{determ-discr-problem}

\subsubsection{Pattern Recognition and Classification}
The problem of image classification into one of two classes $A$ or $B$ can be recast in a form very similar to \eqref{eq:11} and is briefly described here \citep{Novikoff62,Yakubovich63}:
Let $X_k, k=1,2,...N$ denote features, and a corresponding output $y(X_k)$ is of the form
\be y_k=\left\{\begin{array}{ll} 1 & k=1,2,\ldots M\\ -1 & k=M+1,\ldots N\end{array}\right.\label{Frad1}\ee
Suppose that the features are such that these two classes can be separated by a hyperplane in a suitable image of the initial space. That is, the underlying classification model is such that there exist $\theta\in\RR^n$ and $\theta_0\in\RR$ such that 
\be y_k=\theta^T\phi(X_k)+\theta_0\label{fradkov2}\ee
where $\phi(\cdot)$ is a suitable kernel function (regressor) that enables efficient classification. Having exposed the model to a number of  features of known images, the goal is to learn the value of $\theta$ and $\theta_0$ so as to classify any given image into class $A$ or class $B$. That is, the problem is reduced to finding an approximation of the function between $X$ and $y$ based on its values on a finite set. Gradient-type algorithms  in \citep{Novikoff62,Yakubovich63} have been shown to learn the hyperplanes in a finite number of steps \textcolor{blue}{with a prespecified accuracy}.

Yet another approach suggested in \citep{Yakubovich65,Yakubovich66} consists of transforming the above into a dual problem  of finding the intersection of a finite number of half-spaces $$\left\{(\theta,\theta_0): y(X_k) (\theta^T \phi(X_k)+\theta_0)>0\right\}, k=1,2,...N.$$ A gradient-like solution to the above is given by
\begin{equation} \label{f1}
    \theta_{k+1}=\theta_k-\gamma_k y(X_k) \phi(X_k), \theta_{0,k+1}=\theta_{0,k}-\gamma_k y(X_k),
\end{equation}
Several approaches have been proposed to select the size  of the steps (gains) $\gamma_k$. In particular, it is possible to project the current vector of weights onto the boundary hyperplane if the current object is classified incorrectly and take $\gamma_k=0$ otherwise (see details in \citep{Yakubovich66,Bondarko92}).

An alternative approach is to choose a vector of weights $\theta^*$ in such a way that the corresponding hyperplane $\{x: \theta^{*T}\phi(x)+\theta_0=0\}$ is a supporting hyperplane to the convex hull of the available set of vectors ${y(X_k) \phi(X_k)}, k=1,2,...N$, so that the minimum distance from it to the convex hull of classes is maximal. This idea, pioneered in 1964 formed the basis of the celebrated {\it support vector machine (SVM)} method \citep{Vapnik64,Chervonenkis13}.  In the same year, another simple recursive algorithm was also proposed that converges to an optimal supporting hyperplane \citep{Kozinets64} thereby reducing the memory complexity significantly. A min-max based method, {\it MDM} was developed in \cite{Mitchell74} for this problem as well.

A particularly useful method was developed during the 1960s by Bregman \citep, which has become quite popular in recent years in convex optimization and in Machine Learning \citep{Wilson_2016}. In   \cite{Bregman67} (which has more than 1500 citations in Scopus in 2021), Bregman proposed a highly useful notion for a strictly convex function $f(x)$ that later has come to be known as \textit{Bregman divergence}. Currently it is not only of use in convex optimization and associated problems in ML but also in adaptive control. The main idea here is the use of an underlying function $D_{f}(x,y)=f(x)-f(y)-\nabla f(y)^{T}(x-y)$ \bluec{which becomes nonnegative for any $x$ and $y$ if $f$ is convex.} The function $D_f$ is often used either for establishing convergence  or as a Lyapunov function candidate, as discussed in Section \ref{ss:nlp}.  Bregman divergence has also found widespread application in mirror-descent methods in ML. 

It is interesting that a similar problem was addressed in \cite{Gubin67} using a different approach, leading to a number of results on strong convergence and convergence rate in Hilbert space. Additionally, algorithms with incomplete relaxation were proposed and convergence in a finite number of steps was established as well as some applications to Chebyshev approximation and optimal control. This approach has been used in ML as well.

\subsubsection{From pattern recognition to adaptive optimal control}
The link between the above learning methods and adaptive control methods lies in the determination of suitable recursive algorithms so as to 
minimize an underlying loss function as closely as possible.  The papers \cite{Yakubovich65,Yakubovich66} examined this link through ``the method of recurrent goal inequalities" that is based on the reduction of the problem to the solution of a system of inequalities constructed for a given goal function. The proposed gradient-like learning algorithms with deadzone allow one to find a solution to an infinite number of previously not shown inequalities. This in turn allowed Yakubovich to extend his approach to solving adaptive control problems \citep{Yakubovich72,Yakubovich76}, which is stated below.

In \cite{Yakubovich76} an adaptive suboptimal control problem  for a discrete-time linear controlled system  affected by coloured bounded disturbances is studied. Let the controlled system be modeled as follows:
\begin{equation}
  \label{eq:dd1}
  \begin{split}
    y_k = \sum_{i=1}^{n}a_{i}^*y_{k-i} + \sum_{j=1}^{m}b_{j}^*u_{k-j-d}
    + f_k,
      \end{split}
\end{equation}
    \begin{equation}
  \label{eq:dd2}
  \begin{split}
    \sum_{i=0}^{n}c_{i}^*f_{k-i}=g_{k+n},
  \end{split}
\end{equation}
where $y_k,u_k,f_k$ are output, input and disturbances, respectively. Disturbances $f_k$ are generated by the stable filter (\ref{eq:dd2}) with bounded input signal $g_k$ that may take an arbitrary value from the interval $|g_k|\leq G$. Note that equations (\ref{eq:dd1}),(\ref{eq:dd2}) are similar to (\ref{eq:10}) with two differences: first, the filter (\ref{eq:dd2}) is IIR in contrast to the FIR filter in (\ref{eq:10}) and second, disturbances $g_k$ and $f_k$ in  (\ref{eq:dd1}),(\ref{eq:dd2}) are \bluec{deterministic} in contrast to stochastic disturbances $w_k$ in (\ref{eq:10}). 

With the above system \eqref{eq:dd1}-\eqref{eq:dd2}, the following adaptive optimal control problem is now posed. Let the vector of the system parameters $\xi=[a_{i}^*, b_{i}^*,c_{i}^*]$ be unknown while the set $\Xi$ of all admissible values of $\xi$ be known.
Introduce the performance index as follows:
\begin{equation}
  \  \label{eq:dd3}
  \begin{split}
    J_N (u(\cdot)) = &\sup_{k\ge N}\sup_{\xi_f\in\Xi_f}|y_{k}|
  \end{split}
\end{equation}
The problem is to find an admissible adaptive control law not depending on unknown parameters and minimizing the  performance index (\ref{eq:dd3}). It is shown in \cite{Yakubovich76} that if the vector $\xi$ if known, i.e. if $\Xi$   consists of a single element, then the optimal controller would exist and can be represented as
\begin{equation}
  \label{eq:dd4}
  u_k=\theta^{*T}\sigma_k,
  \end{equation}
  where $\theta^*$ is some constant vector explicitly depending on $\xi$ and $\sigma_k$ is the sensor vector. The optimal value (\ref{eq:dd3}) and vector of control parameters (\ref{eq:dd4}) do not depend on $N$ for $N\ge n+r+q$
  
  If the vector $\xi$ is unknown then an adaptation algorithm for the vector $\theta_k$ is proposed allowing to achieve the value of the performance index in the adaptive system arbitrarily close to the optimal value  (\ref{eq:dd3}). It is important that design of the adaptation algorithm is based on the reduction of the problem to the solving an infinite system of the goal inequalities and application of the method \citep{Yakubovich65,Yakubovich66} allowing to solve it in a finite number of steps by a recursive algorithm close to \eqref{f1}. We refer the reader to \citep{Yakubovich76,Fom81,Bondarko92} for further details.

\subsection{Reinforcement Learning/Adaptive Dynamic Programming }


The main problem in RL/ADP can be stated as follows. Let the plant to be controlled be  described by the nonlinear state space equation
\begin{equation}
  \  \label{eq:RL1}
  \begin{split}
    x_{k+1}=f(x_k,u_k)
  \end{split}
\end{equation}
with a control law (also known as a policy)  $\pi=[\mu_1,\mu_2,\ldots]$ as follows
\begin{equation}
  \  \label{eq:RL2}
  \begin{split}
    u_k = \mu_k(x_k).
  \end{split}
\end{equation}
Introduce a performance index in the form of a cost functional as follows:
\begin{equation}
  \  \label{eq:RL3}
  \begin{split}
    J_{\pi}(x_0)  = &\lim_{k\to \infty}\sum_{t=0}^k g(x_t,\mu_t(x_t))
  \end{split}
\end{equation}

The problem is to find $J^*(x)=J_{\pi_*}(x)$, where $J_{\pi_*}(x) = \inf_{\pi}J_{\pi}(x)$.
Based on Bellman's optimality principle, the Value Iteration (VI) process is organized as follows \citep{Bertsekas17}:
\begin{equation}
  \label{eq:RL4}
  \begin{split}
    J_{k+1}(x) = \inf_{u\in U(x)}\left\{g(x,u)+J_k(f(x,u))\right\}.
  \end{split}
\end{equation}
The key step aimed at reducing an overwhelming amount of function evaluations is a neural approximation of the value function as $J_{k}(x) =W_k^{T}\phi(x)$, where $\phi(x)\in R^L$ is truncated set of basis functions, and $W_k$ is a vector of weights that is recursively updated at each $k$.

The difficulty of justifying control based on reinforcement learning under disturbances lies in that the performance index is the average of the integral along the trajectories of the system over the ensemble of disturbances. Averaging requires the use of the Monte Carlo method in one form or another, which inevitably leads to incomplete verification or violation of the stability conditions of the closed loop system. Proving the stability of control systems based on reinforcement learning is a difficult task, and today there are few works where attempts are made to solve it under certain additional assumptions. For example, VI has been thoroughly studied in the setting of Markov Decision Processes with finite states \citep{Wat92,Chang13}. Recent results on   VI for discrete-time dynamical processes evolving in continuous state spaces for nonlinear systems are presented in \citep{Bertsekas17}. Continuous-time counterparts of the   VI related results can  be found  in \cite{BianJiang16,BianJiang21} where the boundedness of all signals and practical stability of the closed loop is established based on neural approximation of both value function and the system Hamiltonian. 
This paper does not address  details of the solutions that have been obtained in this very active and rich research topic.


\section{Solutions}\label{ss:solutions}
This section presents a snapshot of the solutions that have been presented over the last six decades, starting from the results presented in the 70s and 80s. These correspond to stability in continuous-time systems (Section \ref{solutions-stability}), and proceed to learning and parameter estimation (Section \ref{ss:imperfectlearning}). Results proposed for nonlinear systems and nonlinearly parameterized systems follow in sections \ref{ss:nonlinear} and \ref{ss:nlp}, respectively. Similar to deterministic continuous-time systems, the stability results developed in stochastic discrete-time systems are then presented in Section \ref{ss:stochasticstability}. Robustness results obtained starting the 90s are presented in Section IV-F. A cautionary message related to imperfect learning and bursting phenomena is mentioned in Section \ref{ss:bursting}. An inclusion of input constraints in adaptive solutions is discussed in Section \ref{ss:constraints}. Finally, an overview of the assumptions made and the challenges and fundamental tradeoffs encountered in the  evolution of adaptive control are  presented in Section \ref{ss:assn}.
\subsection{Continuous-time Systems-Stability}\label{solutions-stability}

\subsubsection{Algebraic systems}\label{ss:algebraic}
Many problems in adaptive estimation and control may be expressed as 
\be y(t)=\theta^{*T}\phi(t),\label{errormodel1}\ee where $\theta^*,\phi(t)\in\mathbb{R}^N$ represent an unknown parameter and  measurable regressor, respectively, and $y(t)\in\mathbb{R}$ represents a measurable output. This is apparent from \eqref{eq:11} and \eqref{fradkov2}, which corrresponded to estimation and pattern recognition problems in discrete-time. Several examples can be drawn from problems in continuous-time plants as well. One example corresponds to adaptive observers for linear plants, whose output can be represented as an algebraic combination of filtered inputs and outputs (\cite{Narendra2005}, chapter 4). Another example corresponds to a combined-composite approach to adaptive control, discussed at length in papers such as \citep{Slo89,Dua89,Lav09} and more recently a unified approach in \cite{Ort21}, all of which reduce to a plant-model as in \eqref{errormodel1}. 

Given that $\theta^*$ is unknown, we formulate an estimator $\hat y(t)=\theta^T(t)\phi(t)$, where $\hat y(t)\in\mathbb{R}$ is the estimated output and the unknown parameter is estimated as $\theta(t)\in\mathbb{R}^N$. This in turn results two types of errors, a performance error $e_y(t)$ and a learning error $\wt(t)$\footnote{In what follows, we suppress the argument $(t)$ unless needed for emphasis.}
\begin{equation}\label{e:error1}
    e_y=\hat{y}-y,\qquad \wt=\theta-\theta^*
\end{equation}
where the former can be measured but the latter is unknown though adjustable. From  \eqref{errormodel1} and the estimator, it is easy to see that  $e_y$ and $\wt$ are related using a simple regression relation 
\be e_y(t)=\wt^T\phi(t) \label{errormodel1b}.\ee
A common approach for adjusting the estimate $\theta(t)$ at each time $t$ is to determine a rule using all available measurements such that $e_y(t)$ converges towards zero. To do so, a squared loss function
\begin{equation}\label{e:error1_Loss}
    L_1(\theta)=\frac{1}{2}e_y^2
\end{equation} and a corresponding adaptive law for adjusting the parameter error as 
\begin{equation}\label{e:Gradient_Flow_UN}
    \dot{\theta}(t)=-\gamma\nabla_{\theta}L_1(\theta(t)), \qquad \gamma > 0
\end{equation}
is commonly considered \citep{Narendra2005}. It is easy to see that $\nabla_{\theta}L_1(\theta)=\phi e_y$, and therefore \eqref{e:Gradient_Flow_UN} is implementable. Such a gradient-descent approach has found widespread applications in estimation, control, signal processing, and ML. More recently high-order tuners have been proposed for adjusting $\theta$ that uses both gradient and Hessian information
 \citep{Gaudio20AC,Mor2021}.

\subsubsection{Dynamic Systems with States Accessible}
\label{ss:MRAC}
The next class of problems that has been addressed in adaptive control corresponds to plants with all states accessible. We present the solution for the simple case for a scalar input: 
\be \dot x = A_p x + b_pu \label{e:state}\ee where $A_p$ and $b_p$ are unknown, $u$ is the control input and is a scalar, and $x$ is the state and is accessible for measurement. 

\noindent \underline{Matching Condition and Reference Model}
The process of determining an adaptive solution proceeds by first ensuring that the desired solution from the closed-loop system can be described using a reference model. 
For the plant in \eqref{adaptive1}, this reference model takes the form 
\be \dot x_m = A_m x_m+b_m r \label{e:refmodel}\ee 
 and is such that the state $x_m(t)$ encapsulates the desired solution expected from the controlled plant. This can be accomplished by choosing a reference input $r$, $A_m$ to be a Hurwitz matrix, $(A_m,b_m)$ is controllable so that together they produce an $x_m(t)$ that approximates the signal that the plant is required to track. 

With the reference model  chosen as above, the next step pertains to Matching Conditions \citep{Narendra2005}. These ensure that a controller with fixed parameter exists, which guarantees that the closed-loop system matches the reference model. In particular, for the plant in \eqref{e:state}, a control input of the form
\be u(t)=\theta^{*T}x(t)+k^* r(t) \label{e:fixedcontrol}\ee
guarantees this match provided $\theta^*$ and $k^*$ solve the following which is denoted as Matching Conditions:
\bena A_p+b_p\theta^{*T}&=& A_m\nonumber\\
b_p k^*&=& b_m \label{e:match}\eena
This corresponds to the Algebraic Part of the problem described in Section \ref{ss:cep}.

The final step is the analytic part, the rule for estimating the unknown parameters $\theta^*$ and $k^*$ and the corresponding adaptive control input that replaces the input choice in \eqref{e:fixedcontrol}. These solutions are given by
\begin{eqnarray}
u&=&\theta^T(t)x+k(t)r\label{adaptive1}\\
\dot \theta &=& - sign(k^*) \Gamma_\theta (e^TPb_m)x\label{adaptive2}\\
\dot k& =& - sign(k^*) \gamma_k (e^TPb_m)r\label{adaptive3}\eea
where $\Gamma_\theta>0$ is a positive definite matrix, $\gamma_k>0$ is a positive constant, $e=x-x_m$, and $P=P^T\in\mathbb{R}^{n\times n}$ is a positive definite matrix that solves the Lyapunov equation \be A_m^TP+PA_m=-Q\label{Lyap}\ee with a positive definite matrix $Q=Q^T\in\mathbb{R}^{n\times n}$. It can be shown that 
\be V=e^TPe+|k^*|\left[(\theta-\theta^*)^T\Gamma^{-1}(\theta-\theta^*)+(1/\gamma_{k})(k-k^*)^2\right]\label{lyap-state}\ee
is a Lyapunov function with $\dot V=-e^TQe$ and that $\lm{t}e(t)=0$. The reader is referred to Chapter 3 in  \cite{Narendra2005} for further details.

It should be noted that the adaptation rules in \eqref{adaptive2}-\eqref{adaptive3} can also be expressed as the gradient of a loss function \citep{Gaudio20AC}
\begin{equation}\label{e:error2_Loss}
    L_2(\bar\theta)=\frac{d}{dt}\left\{\frac{e^TPe}{2}\right\}+\frac{e^TQe}{2},
\end{equation}
where $\bar\theta=[\theta^T,k]^T$, and it is assumed that $k^*>0$ for ease of exposition.
  It is noted that this loss function $L_2$ differs from that in \eqref{e:error1_Loss}, and includes an additional component that reflects the dynamics in the system. 
  It is easy to see that 
\begin{equation}\label{e:Gradient_Flow2}
    \dot{\bar\theta}(t)=-\Gamma\nabla_{\bar\theta}L_2(\theta(t)), \qquad \Gamma > 0,
\end{equation}and is implementable as  $\nabla_{\theta}L_2(\theta)=\phi e^TPb_m$, can be computed at each time $t$, where $\phi=[x_p^T, r]^T$. 

The matching condition  \eqref{e:match} is akin to the controllability condition, albeit somewhat stronger, as it requires the existence of a $\theta^*$ for a known Hurwitz matrix $A_m$ \citep{Lavretsky2013,Narendra2005}. The other requirement is that the sign of $k^*$ needs to be known, which is required to ensure that $V$ is a Lyapunov function. 
\subsubsection{Adaptive Observers}
The adaptive control solution in \eqref{adaptive1}-\eqref{adaptive2} in section \ref{ss:MRAC} required that the state $x(t)$ be available for measurement at each $t$. A central challenge in developing adaptive solutions for plants whose states are not accessible is the simultaneous generation of estimates of both states and parameters in real-time. Unlike the Kalman Filter in the stochastic case or the Luenberger observer in the deterministic case, the problem becomes significantly more complex, as state estimates require plant parameters and parameter estimation is facilitated when states are accessible. This loop is broken using a non-minimal representation of the plant, leading to a tractable observer design. Starting with a plant model as in \eqref{plant-linear}, a state-representation of the same can be derived as is given by \cite{Luders_1974} 
\bea \dot\omega_1&=&\Lambda \omega_1+\ell u\nonumber\\
\dot\omega_2&=& \Lambda \omega_2+\ell y\label{e:adaobs}\\
y&=&\theta_1^T\omega_1+\theta_2^T\omega_2\nonumber\eea
where $\omega=[\omega_1^T,\omega_2^T]^T$ is a nonminimal state of the plant transfer function $W_p(s)$ between the input $u$ and the output $y$. $\Lambda\in\RR^{n\times n}$ is a Hurwitz matrix and $(\Lambda,\ell)$ is controllable and are known parameters. Assuming that $W_p(s)$ has $n$ poles and $m$ coprime zeros, in contrast to a minimal  $n^{th}$-order representation, Eq. \eqref{e:adaobs} is nonminimal and has $2n$ states. The adaptive observer leverages Eq. \eqref{e:adaobs} and generates a state estimate $\widehat\omega$ and a plant estimate $\wh$ as follows:
\bea \dot{\hat\omega}_1&=&\Lambda \hat\omega_1+\ell u\nonumber\\
\dot{\hat\omega}_2&=& \Lambda \hat\omega_2+\ell y\label{e:adaobs2}\\
\hat y&=&\wh_1^T\hat\omega_1+\wh_2^T\hat\omega_2\nonumber\eea
where $\wh=[\wh_1^T,\wh_2^T]^T$ and $\hat\omega=[\hat\omega_1^T,\hat\omega_2^T]^T$. The adaptive law that adjusts the parameter estimates is chosen as 
\be \dot\wh=-\Gamma \left(\hat y_p-y_p\right)\widehat\omega\label{e:adlawobs}\ee
where $\Gamma$ is a known symmetric, positive definite matrix. 

Analytical guarantees of stability of the parameter estimate $\wh$ in \eqref{e:adaobs2} and \eqref{e:adlawobs} as well as asymptotic convergence of $\wh(t)$ to $\theta$ can be found in \cite{Morgan_1977,Narendra_1987}. Necessary and sufficient conditions for this convergence requires that the regressor $\hat\omega$ be persistently exciting. Several results also exist in ensuring accelerated convergence of these estimates \citep{Lion_1967,Kreisselmeier_1977,Jenkins_2019,Aranovskiy_2019,Ort21,Gaudio_2020}) using  matrix regressors, a time-varying learning rate for $\Gamma$, and dynamic regressor extension and mixing.

\subsubsection{Adaptive Controllers with Output Feedback - A special case}
The two assumptions made in the development of adaptive systems in Section \ref{ss:MRAC} include matching conditions and the availability of states of the underlying dynamic system at each instant $t$. Both are often violated in many problems, which led to the development of adaptive systems when only partial measurements are available. With the focus primarily on linear time-invariant (LTI) plants, the first challenge was to address the problem of separation principle employed in control of LTI plants \citep{Kailath1980,CTChen1984}. The idea therein is to allow a simultaneous estimation of states using an observer and a feedback control using state estimates with a linear quadratic regulator to be implemented and allow them both to proceed simultaneously in real-time and guarantee stability of the closed-loop system. The challenge in the current context is that parameters are unknown, introducing an additional estimate, of the plant parameter, to be generated in real-time. In contrast to the classical problem where the closed-loop remains linear, the simultaneous problem of generating the parameter estimate to determine the controller and the feedback control to ensure the generation of well-behaved parameter estimates introduced intractable challenges.

The starting point is an input-output representation of the plant model as in \eqref{plant-linear}. Recognizing that estimation and control are duals of each other \citep{Fel60}, a similar nonminimal representation of the plant as in \eqref{e:adaobs} was used as the starting point to decouple the estimation of the state from the design of the control input. In particular, an adaptive control input of the form 
\be u(t)=\theta_c^T(t)\omega(t)+k(t)r(t)\label{e:adcon}\ee
enabled a tractable problem formulation, where $\omega(t)$ is generated as in \eqref{e:adaobs}. 
The added advantage of the nonminimal representation is that it ensures the existence of a solution that matches the controlled plant using \eqref{e:adcon} to that of the reference model. That is, the existence of a control parameter $\theta^*$ and $k^*$ such that 
\be u(t)=\theta^{*T}\omega(t)+k^*r(t)\label{e:admatch}\ee
ensured that the closed-loop transfer function from $r$ to $y$ matched that of a reference model with a transfer function $W_m(s)$, specified as
\be y_m(t)=W_m(s)r(t)\label{output-refmodel}\ee 
That is, the controller in \eqref{e:admatch} is guaranteed to exist for which the output error $e_y=y_p-y_m$ has a limiting property of $\lm{t}e_y(t)=0$. For this purpose, the well known Bezout Identity \citep{Kailath1980} and the requirement that $W_p(s)$ has stable zeros was leveraged. 

When the adaptive controller as in \eqref{e:adcon} is used, the plant model in \eqref{plant-linear} and the existence of $\theta^*$ and $k^*$ that guarantee that the output error $e_y(t)$ goes to zero, leads to an error model of the form 
\be e_y=(1/k^*)W_m(s) [\bar\wt^T\phi] \label{e:SPR}\ee
where $\phi=[\omega^T, r]^T$, $\bar\wt=[(\theta-\theta^*)^T, (k-k^*)]^T$.

The problem of determining the adaptive rule for adjusting $\bar\wt$ was solved in a very elegant manner when the relative degree, i.e. the net-order of $W_m(s)=1$. In this case,  a fundamental systems concept of strictly positive real (SPR) transfer function as well as an elegant tool known as  Kalman-Yakubovich Lemma (KYL) \citep{Yakubovich1962,Kalman1963,Mey65,Lef63,Nar73,Anderson_1982} can be leveraged. 
This KYL was first
proposed  by Yakubovich \citep{Yakubovich1962} and extended by Kalman \citep{Kalman1963}, which came out of stability theory of nonlinear systems, Popov's absolute stability, and the Circle Criterion \citep{Nar73}. This is briefly described below.



\noindent\underline{Strictly Positive Real Functions:}
The concept of positive realness arose in the context of stability of a class of linear systems with an algebraic nonlinearity in feedback. It was shown, notably by Popov, that under certain conditions on the frequency response of the linear system, that a Lyapunov function can be shown to exist. The KYL establishes the relation between these frequency domain conditions and the existence of the Lyapunov function. Both the definition of rational SPR functions and the KYL are listed in the appendix.

Using the KYL, the following adaptive laws are proposed for the adjustment of the control parameters:
\bea
\dot\theta& =& -sign(k^*)e_y\omega\label{ad20}\\
\dot k &=& -sign(k^*)e_y r\label{ad21}\eea
It can be shown that \be V=e^TPe+(1/|k^*|)\left(||(\theta_c-\theta^*)||^2+|(k-k^*)|^2\right)\label{e:LyapSPR}\ee 
is a Lyapunov function where $P$ is the solution of the KYL for the realization $\{A_m,b,c\}$ of $W_m(s)$ which is SPR. This follows by first noting that
\be 
\dot V = -e^T(A_m^TP+PA_m)e + 2e^TPb(\wt^T\omega+\tilde k r)-\dot{\wt}^Te_y\omega-\dot{\tilde k}e_y r.\nonumber\ee
Since $W_m(s)$ is SPR, the use of the KYL applied to $W_m(s)$ together with the adaptive laws in \eqref{ad20}-\eqref{ad21} and \eqref{e:LyapSPR} causes the second term to cancel out the 3rd and 4th terms and hence that
  $\dot V=-e^TQe\leq 0$. The structure of the adaptive controller in \eqref{e:adcon} guarantees that $e_y$, $\theta_c$, $k$, $\omega$, $y_p$, and $u$ are bounded and that $\lm{t}e_y(t)=0$. Additions of positive definite gains to \eqref{ad20} and \eqref{ad21} as in \eqref{adaptive2}-\eqref{adaptive3} are straight forward.
 
The choice of the adaptive laws as in \eqref{ad20}-\eqref{ad21} centrally depended on the KYL which in turn required that $W_m(s)$ be SPR. A SPR transfer function \citep{Narendra2005} leads to the requirement that the relative degree, the difference between the number of poles and zeros of $W_p(s)$, is unity, and has stable zeros (zeros only in $Re[s]<0$), also defined as \textit{hyperminimum-phase} \citep{Fradkov74}. Qualitatively, it implies that a stable adjustment rule for the parameter should be based on loss functions that does not significantly lag the times at which new data comes into the system. For a general case when the relative degree of $W_p(s)$ exceeded unity, it posed a significant stability problem, as it was clear that the same simple adaptive laws as in \eqref{ad20}-\eqref{ad21} will no longer suffice as the corresponding transfer function $W_m(s)$ of the reference model cannot be made SPR.

A final note about the assumptions made about the plant-model in \eqref{plant-linear} is in order. For the controller in \eqref{e:admatch} to allow the closed-loop system to match the reference model in \eqref{output-refmodel} for any reference input $r(t)$, a reference model $W_m(s)$ with the same order and net-order as that of $W_p(s)$ needs to be chosen, which implies that the order and net-order of the plant need to be known. Determination of a Lyapunov function requires that the sign of $k^*$ be known. Finally, the model-matching starting with a non-minimal representation of the plant required stable pole-zero cancellations, which necessitated the zeros to be stable.

\subsubsection{Adaptive Controllers with Output Feedback - Passification approach}
In some cases the structure of the adaptive controller may be significantly  simplified avoiding usage of reference model or adaptive observer. This approach  is based on the so called {\it passification lemma} - a feedback version of KYL \citep{Fradkov74,Fradkov03}. For simplicity consider  the case of stabilization  $r(t)=0$. Let vector $B\in R^n$, $n\times l$ matrix $C$ and vector $g\in R^l$ be given.

\noindent {\it Lemma (Passification)}: Consider matrix relations 
\bea A_{\theta}^TP+PA_{\theta}&<& 0, A_{\theta}= A+B\theta C\nonumber\\ PB=Cg\label{e:fkyl}\eea
There exist a symmetric positive definite matrix $P$ and a vector $\theta$ satisfying (\ref{e:fkyl}) if and only if the transfer function $Z_g(s)=(Cg)^T(sI-A)^{-1}B$ is {\it hyperminimum-phase}.

Based on the above lemma it can be proven that the control plant 
\bea \dot x = Ax+Bu, y = C^Tx \label{e:fkyl2}\eea
can be stabilized by adaptive controller
\bea u = \theta^Ty, \dot\theta=-\Gamma (g^Ty)y \label{e:fkyl3}\eea
in the sense that $x(t)\to 0, \theta(t)\to const$ as $t\to\infty$, if
$Z_g(s)$ is hyperminimum-phase. Moreover the property of hyperminimum-phaseness is necessary and sufficient for existence of Lyapunov function
\bea V(x,\theta) = x^TPx + 0.5(\theta-\theta^*)^T\Gamma^{-1}(\theta-\theta^*)
 \label{e:fkyl4}\eea
such that $V(x,\theta)>0$ for $x\neq 0, \theta \neq \theta^*$ and $\dot V<0$ for $x\neq 0$.
Extensions and applications of passification approach can be found in \citep{AndrFr06, AndrSeliv18}.

\subsubsection{Adaptive Controllers with Output Feedback - The General case}\label{s:general}
Extensions to a general case with output feedback have been proposed using several novel tools including an augmented error approach \citep{Narendra_1989}, backstepping \citep{Krstic_1995}, averaging theory \citep{And86}, and high-order tuners \citep{Annaswamy2003}. In all cases, the complexity of the adaptive controller is increased, as the error model in \eqref{e:SPR} does not permit the realizations of simple loss functions as in $L_i(\theta)$, $i=1,2$. 

 Over the years, several solutions have been proposed in the literature to address this problem, of which two are briefly summarized below. In all cases, the zeros of $W_p(s)$ are required to be stable.

\noindent\underline{Augmented Error Approach:}
The problem is to convert an error model that is of the form
\be e_1=W_m(s)[\tilde\theta^T\omega]\nonumber\ee
where $W_m(s)$ is not SPR to one where the transfer function between the parameter error and the output error is SPR. Towards this end, an auxiliary error $e_2$ is added to $e_1$, where
\be e_2=[\theta^TW_m-W_m\theta^T]\omega \label{e:aux}\ee
where $\theta$ is the adjustable parameter and $\wt$ is the corresponding parameter error. It can be shown that the resulting augmented error $\epsilon_1=e_1+e_2$ has a simple error model structure of the form
\be \epsilon_1=\wt^T\zeta, \qquad \zeta=W_m(s)\omega\label{e:augment}\ee
This in turn allows a simple adjustment rule 
\be \dot{\wt}=-m(t)\epsilon\zeta\label{adlaw-general}\ee 
where $m(t)$ is a suitably chosen normalizing signal that guarantees not only that $\theta(t)$ is bounded but also that $\dot\theta\in{\cal L}_2$. These two properties of the parameter are then suitably leveraged to guarantee that the closed loop system with the controller defined by \eqref{e:adcon} and the adaptive law specified by \eqref{e:aux}-\eqref{adlaw-general} has bounded solutions and that the output error $e_y(t)$ converges to zero asymptotically. BIBO properties of linear systems, almost time-invariant systems, and minimum-phase systems as well as order arguments are leveraged in this proof of stability \citep{Narendra1980,Narendra2005}. As in the case of the special case when the net-order of $W_p(s)=1$, the requirements that the order, net-order, the sign of the high-frequency gain be known, and that the zeros of $W_p(s)$ are stable are all needed in this case as well.

\noindent\underline{High-order Tuners:}
The starting point is a plant model of the form
\begin{equation}
	y(t) = W_p(s)[u(t)],\; W_p(s) = \frac{Z_p(s)}{R_p(s)}\label{eq:0}
\end{equation}
with a relative degree $m$ that is known, an order $n$ that is unknown, with all zeros in $C^-$, and a high frequency gain unity. It is well known that this plant can be stabilized using
\begin{equation}
	u=-\frac{p(s)}{(s+z_c)^{m-1}}u - k_1 y
	\label{eq:1}
\end{equation}
where
\begin{equation}
	p(s) = k_{21} + k_{22}s + \cdots + k_{2(m-1)}s^{m-2}
	\label{eq:2}
\end{equation}
for suitable values of $k_c$ and $k_{2j}$ in \eqref{eq:2}.  This follows from the fact that the closed-loop transfer function is of the form
\begin{equation}
	W_{cl}(s) = \frac{(s+z_c)^{m-1}Z_p(s)}{R_p(s)p_c(s) + k_1(s+z_c)^{m-1}Z_p(s)}
\end{equation}
where
\begin{equation}
	p_c(s) = (s + z_c)^{m-1} + p(s).
\end{equation}
For a large $k_1$, the $n+m-1$ poles of $W_{cl}(s)$ become close to the zeros of $(s+z_c)^{m-1}Z_p(s)$ and other $m$ stable locations for suitable values of $k_{2i}, i=1,\ldots, m-1$ \citep{Annaswamy2003}. We now utilize the controller structure as in \eqref{eq:1}-\eqref{eq:2} to describe an alternate type of adaptive controller. 

A time-domain representation of \eqref{eq:1}-\eqref{eq:2} is given by
\begin{equation}
	\dot\omega_1 = \Lambda\omega_1 + \ell u
\end{equation}
\begin{equation}
	u = - k_2^T\omega_1 - k_1 y + r
	\label{eq:3}
\end{equation}
where $\Lambda\in\RR^{m\times m}$, $(\Lambda, \ell)$ is controllable, and
\begin{equation}
	k_2^T(sI - \Lambda)^{-1}\ell = \frac{p(s)}{(s+z_c)^{m-1}}.
\end{equation}
As $k_1$ and $k_2$ are unknown when the parameters of $W_p(s)$ are unknown, an adaptive controller corresponding to \eqref{eq:3} is given by
\begin{equation}
	u=k_2^T(t)\omega_1 + k_1(t)y + r.
\end{equation}
Expressing the control parameters as $k_1(t) = k_1^* + \tilde k_1(t)$, $k_2(t) = k_2^* + \tilde{k}_2(t)$, $\omega = [\omega_1^T, y]^T$, $\tilde{k} = [\tilde k_1^T, \tilde k_2^T]^T$, the closed-loop system equations can be described as
\begin{equation}
	y = W_{cl}(s)(\tilde k^T \omega) + r.
\end{equation}
$W_{cl}(s)$ is not strictly positive real (SPR), but has stable poles, stable zeros, and is of relative degree $m$. Due to these properties, it is reasonable to assume that one can find a strictly positive real transfer function of the form
\begin{equation}
	W_m(s) = W_{cl}(s) (s+a)^{m-1}.
\end{equation}
To enable the realization of $W_m(s)$ in closed-loop, we choose the control input, instead of $k^T(t)\omega(t) + r(t)$, as follows:
\begin{align}
	u(t) &= (s+a)^{m-1}[k^T\omega'(t) + r']\label{eq:ht4}\\
	\omega'(t) &= \frac{1}{(s+a)^{m-1}}[\omega(t)]\label{eq:ht5}\\
	r'(t) &= \frac{1}{(s+a)^{m-1}}[r(t)]
\end{align}
This will lead to
\begin{equation}
	y = W_m(s)(\tilde k^T\omega' + r').
	\label{eq:6}
\end{equation}
Now, the problem is to realize \eqref{eq:ht4} without explicitly differentiating any signal. Let $p = m - 1$. Using binomial expansion and the chain rule for differentiation, we obtain that
\bena
	u &=& k^T d_0 + p\dot k^T d_1 + \cdots + (pc_i){k^{(i)}}^T d_i + \cdots + 
	\\&& \hbox{}\;\;+ {k^{(p)}}^T d_p + r,\\
	d_i(t)& =& \left[\frac{1}{(s+a)^i}\right][\omega(t)], \; i = 1, \ldots, p
\eena
Note that all terms involving $k$ and $d_i$ are realizable. So, the only remaining piece is the realization of derivatives of $k$ to $p$th order.

The overall problem can be summarized as follows: Given the closed-loop system in Eq. \eqref{eq:6} where $\omega'$ is given by \eqref{eq:ht5}, determine an adaptive law for adjusting $k$ so that it is differentiable $p$ times and all the signals in the loop are bounded. The time-domain representation of the error model in \eqref{eq:6} is given by
\begin{equation}
	\dot e = A_s e + b_s(k - k^*)^T\omega', \; e_1 = h_s^T e,
\end{equation}
where
\begin{equation*}
	h_s^T(sI - A_s)^{-1}b_s = W_m(s).
\end{equation*}
Since $W_m(s)$ is SPR, we have that
\begin{equation}
	A_s^T P_s + P_s A_s = -Q\leq 0,\; P_sb_s = h_s
\end{equation}
We note that $\omega'$ is differentiable $p$ times. In what follows, $\omega_i'$ and $k_i$ denote the $i$th element of a vector $\omega'$ and $k$, respectively.

Using the high-order tuners developed in \cite{Morse_1992}, the following adaptive law is suggested for adjusting $k$,
\begin{align}
	\dot k' &= - e_1\omega'\label{eq:7}\\
	\dot x_i &= (Ax_i + bk_i')f(\omega_i'), \; f(x) = 1 + \mu x^2\\
	k_i &= c^T x_i\\
	c^T(sI - A)^{-1}b &= \frac{\alpha(0)}{\alpha(s)}\label{eq:8}
\end{align}
and $\alpha(s)$ is an arbitrary stable polynomial of degree $p$. The choice of $k$ as in Eqs. \eqref{eq:7}-\eqref{eq:8} guarantees that $k$ is differentiable $p$ times.

\subsection{Learning and Persistent Excitation}\label{ss:imperfectlearning}
The focus of all problems addressed in Sections \ref{solutions-stability} is to bring the performance error $e_y(t)$ or $e(t)$ to zero. This performance corresponds to either successful  output estimation or tracking, both of which are reflected in the choice of the underlying loss function. However an additional goal in many adaptive systems is to learn the underlying parameters. As is clear from all preceding discussions, the hallmark of all adaptive control problems is the inclusion of a parameter estimation algorithm. In addition to ensuring that the closed-loop system is bounded and that the performance errors are brought to zero, all adaptive systems attempt to learn the underlying parameters, with the goal that the parameter error $\theta-\theta^*$ is reduced if not brought to zero. We discuss two important aspects under which this learning, i.e. reduction of parameter error to zero, occurs.

The first is the necessary and sufficient condition under which learning occurs:
\begin{definition}[\hspace{1sp}\cite{Narendra2005}]\label{d:PE}
    A bounded function $\phi:[t_0,\infty)\rightarrow\mathbb{R}^N$ is persistently exciting (PE) if there exists $T\hspace{-.045cm}>\hspace{-.045cm}0$ and $\alpha\hspace{-.045cm}>\hspace{-.045cm}0$ such that
    \begin{equation*}
        \int_t^{t+T}\phi(\tau)\phi^T(\tau)d\tau\geq\alpha I,\quad \forall t\geq t_0.
    \end{equation*}
\end{definition}
It has been shown in \cite{Morgan_1977,Narendra2005} that this leads to convergence of the parameter error in algebraic systems, dynamic systems with states accessible, and in those with output feedback. Several books and papers have delineated properties of the exogenous signals in a control system that ensures the underlying regressor $\phi$ is persistently exciting \citep{Narendra_1987,Narendra_1989,Sastry_1989,Boyd_1983}. It should be noted that this property creates a rank $N$ matrix over an interval despite the fact that the integrand is of rank one at any instant $\tau$. Conditions that ensure parameter learning with high-order tuners in Eqs. \eqref{eq:7}-\eqref{eq:8} are established in \citep{ortega1993}. This necessary and sufficient condition on the underlying regressor is shown to lead to several desirable properties of the adaptive system, including lack of bursting \citep{And85,Mor77,For81,Nar87} and uniform asymptotic stability and robustness to disturbances \citep{Nar86,Jenkins2018}.

The second is the important observation that persistent excitation is not required for satisfactory performance of the adaptive system; both output estimation and tracking, typical goals in system estimation and control, can be achieved without relying on learning. That is, a guaranteed safe behavior of the controlled system can be assured in real-time prior to reaching the learning goal. This guarantee in the presence of imperfect learning is essential, and suggests that for real-time decision making, {\em control for learning} is the practical goal in contrast to {learning for control}.

\subsection{Nonlinear Systems}\label{ss:nonlinear}

All of the discussions above pertain to the linear plant model in \eqref{plant-linear}. We now return to the original problem in \eqref{plant1}, where we assume that the unknown parameter 
$\theta \in \Xi \subset R^N$ is a vector of unknown parameters belonging to an a priori known set  $\Xi$. Let the control goal be
\begin{equation} \label{F5-2}
Q_t\le \Delta ~for ~t\ge t_* 
\end{equation}
where $Q_t=Q[x(s),u(s); 0\le s\le t]$  is the objective functional.
The task is: to find a two-level control law
\begin{equation}\label{F5-3}
\begin{array}{c}
u(t) = U_t[y(s), u(s), \theta(s); 0 \le s < t],\\
\theta(t) =\Theta_t(y(s), u(s), \theta(s); 0 \le s < t)
\end{array}\end{equation}
such that in the closed loop control system \eqref{plant1} and \eqref{F5-3} meets the goal  \eqref{F5-2} and
its trajectories remain in the sets $x(t)\in D_x, u(t)\in D_u, \xi(t)\in D_{\xi}$
for any $\xi \in \Xi$ and $(x(0),\theta(0))\in Q_0$ where
$Q_0\subset D_x\times D_{\theta}$ is a prespecified set. Here $\theta(t)$ is a vector of
adjustable parameters. Note that operators $Q_t, U_t$ and $\Theta_t$ are all nonanticipative.

A number of studies were aimed at relaxation of matching conditions for nonlinearities. 
A breakthrough was made in the end of 1980s by several groups. Further development made in \cite{Krstic95} lead to an elegant technology of iterative control design called ``backstepping design''. The number of the papers using it for adaptive control was growing rapidly and exceeded one thousand during the decade 2011-2020 with about a quarter dedicated to nonlinear adaptive control
(counting by the number of papers in the Web of Science database with the terms ``backstepping  AND adaptive  AND nonlinear'' in the paper title).
An approach advocated in \cite{Krstic95} and related papers is based on  the application of the backstepping procedure directly to model \eqref{plant1} expressed in a standardized canonical form. An alternative approach, proposed 
in \cite{Marino91,Marino93}, assumes the use of special filters, which are part of the adaptive controller, that make 
it possible to transform the model in \eqref{plant1} to an ``adaptive observer canonical form'', and then apply the backstepping procedure to the transformed plant model. The class of problems considered in the above papers was expanded further in \cite{Seto_1994} to include triangular structures.

A number of approaches to adaptive control of nonlinear systems are based on approximation of 
 nonlinear right hand sides by linear ones. There are only a few publications with explicit
formulations of dynamic properties of the overall system, e.g. the paper  \cite{WenHill90}, where
reduction of the nonlinear model is made by standard linearization via finite differences;
There are a few results dealing with high gain linear controllers for nonlinear systems
\citep{Gusev88,Marino85}. 


Finally, it should be pointed out that adaptive control of nonlinear systems have also employed fundamental tools  such as absolute stability \citep{Haddad01,FrLipkovich15}, passivity \citep{Ast07}  and passification \citep{Fom81,SeronHillFradkov95,AndrSeliv20}. Also noteworthy is a related general approach based on immersion and invariance \citep{Ast03}.

\subsubsection{Nonlinear Control with Neural Networks}
Since the late 1980s, there has been a rapid growth in the number of works devoted to the adaptive control of nonlinear systems based on learning and neural networks. The basic principles of using artificial neural networks in control problems were formulated in the seminal article \citep{Narendra_1990}, which received more than 5000 citations over two decades. Neural networks are widely used as a means of approximating nonlinear functions for learning and control by many ways. In a number of works, neural networks are used to approximate the right-hand sides of the system. For example, in the work \citep{Polycarpou96} 
 which got more than one thousand citations it is proposed   to approximate a nonlinear scalar function $f(x))$ in the 2nd order equation
$\dot x_1=x_2+f(x_1), \dot x_2=u$ by a linear combination of radial basis neural network functions with tunable weights and then to adjust the weights based on the online measurements. Adaptation algorithm for weights is based on the Lyapunov function which is quadratic with respect to both plant state and tunable weights. An extended adaptation algorithm for a class of $nth$ order nonlinear systems was proposed in \cite{Sanner_1992} and a more general case was studied in \cite{Lewis96}.

In a number of works  deep (multilayer) neural networks were used to approximate the right-hand side of \eqref{F5-3}, but the only weights that are adjusted are in the outer layer of the network, with the weights in all remaining layer fixed \citep{Hov08,Ren10,Rov94,Lavretsky2013}. This makes the problem tractable as the underlying Lyapunov function can still be chosen to be quadratic. Very few solutions have been provided when the hidden layers are also adjusted, as it makes it very difficult to prove convergence \citep{Lewis96,Patkar2020}. Another approach that has been used is an approximation of Lyapunov functions in closed-loop systems using Neural Networks together with its derivative along the system, that must satisfy the inequalities justifying the stability of the closed-loop \citep{Chang20,Yu_1998}. However, the verification of the fulfillment of inequalities should be carried out in the whole space or in a representative set as even a small violation of the inequality $\dot V<0$  may lead to an incomplete verification or violation of the stability conditions.


\subsubsection{Parameter Learning in Nonlinear Systems}
Conditions for parameter learning have been investigated at length in nonlinear systems as well, by posing the underlying problem as the uniform asymptotic stability (UAS) of a nonlinear differential equation
\be \dot x = F(x,t) \label{e:gas}\ee 
where $x$ corresponds to the underlying parameter error. We refer the reader to \cite{loria2005} for an excellent exposition of the underlying results as well as the references therein for details of this topic. We briefly summarize the idea below: 
The challenge in all adaptive systems including those that arise in the context of control of linear plants is that the underlying Lyapunov function can only be shown to be negative semi-definite, while parameter convergence, i.e. UAS of \eqref{e:gas} requires negative definiteness of a Lyapunov function. This challenge is tackled in \cite{Mor77} by applying uniform observability properties of linear systems. In \cite{loria2005}, a new definition of persistent excitation and the use of Matrosov’s theorem are utilized to achieve UAS for nonlinear systems. Matrosov’s theorem can be viewed as an invariance principle for nonautonomous systems, and revolves around constructing an auxiliary function on top of a Lyapunov function, with a nonzero derivative on the set where the Lyapunov’s function has a derivative that is zero. These tools are shown to be applicable for a class of nonholonomic systems.

\subsection{Nonlinearly Parameterized Systems}\label{ss:nlp}
All of the problems described thus far, both in deterministic and in stochastic systems have assumed that the parametric uncertainties appear linearly. 
 A class of problems that have relaxed this assumptions can be found in \citep{Fradkov79,AndrFr21,Ort96,Annaswamy1998,Annaswamy98,Loh_1999,Frad01} and have provided solutions for problems when parameters occur nonlinearly. The starting point for these solutions is \textit{speed-gradient method} \citep{Fradkov79,AndrFr21} which not only works for nonlinear systems such as in \eqref{plant1} but also for nonlinearly parameterized systems. It is assumed that a parametric stabilizing feedback law $u=U_*(x,\theta,t)$
 is known such that if  $\theta=\theta_*$ then $Q(x,t)$ along trajectories of 
the closed loop (\ref{plant1}) are such that $w(x,\theta,t,\xi)=\partial Q/\partial t+\partial Q/\partial x F(x,U_*(x,\theta,t,\xi)$ is negative 
definite in $x$. Then the speed-gradient control is designed as follows:
\begin{equation} \label{F5-5}
\dot\theta=-\Gamma\nabla_{\theta} w(x,\theta,t),
\end{equation}
where $\Gamma=\Gamma^{T}>0$ is positive definite matrix gain with a Lyapunov function
\begin{equation} \label{F5-1b}
V(x,\theta,t)=Q(x,t)+(\theta-\theta_*)^T\Gamma^{-1}(\theta-\theta_*).
\end{equation}
An assumption that $w(x,\theta,t)$ is convex in $\theta$  is required. Inspired by this approach, further extensions were reported in \citep{Ort96,Annaswamy1998,Loh_1999} and is briefly summarized below.

Suppose the underlying nonlinear system is of the form
\be \label{plant}
    \dot X_p = A_pX_p+b(f(X_p,\theta) X_p + u)
\ee
where $X_p \in \mathbb{R}^n$ is the plant state assumed accessible for measurement,  $A_p\in\RR^{n\times n}$, $b \in \mathbb{R}^{n}$, $u \in \mathbb{R}^{m}$ is the control input. The function $f$ is nonlinear not only with respect to $X_p$ but also with respect to the parameter $\theta$. Typical examples of such nonlinear parameterizations are all types of neural networks including deep networks and radial basis functions, and all physical systems with complex constitutive relations \citep{Annaswamy98}. The main difficulty posed by the nonlinearity in $\theta$ is briefly explained below.

The structure of the plant dynamics in \eqref{plant} suggests that when $A_p$ and $\theta$ are known, a control input of the form 
\be 
u = -f(X_p,\theta)+\alpha^T X_p+r \label{nparam-input}\ee
leads to a closed-loop system with BIBO properties, given by
\ben \dot X_p=A_mX_p+b_r\een 
where $\alpha$ satisfies the matching condition
$A_p+b_p\alpha^T=A_m$, and $A_m$ is a Hurwitz matrix. Inspired by the control structure in \eqref{nparam-input}, the adaptive counterpart of the same that attempts to control \eqref{plant} and learn the parameters $\alpha$ and $\theta$ through an estimation process is given by 
\be 
u = -f(X_p,\hat\theta)+\hat\alpha^T X_p+r \label{nparam-inputa}\ee
One can now derive an error equation 
\be 
\dot e = A_m e + b(f-\hat f+\hat\alpha^TX_p)\label{np-error}\ee
where $e=X_p-X_m$, $\hat f=f(X_p,\hat\theta)$, and $X_m$ is the state of a reference model
\be \dot X_m=A_mX_m+br \label{refmodel}\ee
With a few transformations, the vector equation in \eqref{np-error} can be reduced to a scalar error equation 
\be \dot e_c = -k e_c +f-\hat f \label{np-error2}\ee
where $e_c=h^Te$ and $k>0$ \citep{Annaswamy1998}. If one were to choose a standard quadratic Lyapunov function candidate $V=e_c^2+\wt^2$, where $\wt=\hat\theta-\theta$, its time-derivative is of the form 
\be \dot V = -ke_c^2+e_c\left[f-\hat f+\wt\nabla f_{\hat\theta}\right]\label{dotV-np}\ee
Unlike the case when $f$ is linear in $\theta$, where one could choose an adaptive law for adjusting $\wt$ so that the term within the brackets will become identically zero, one cannot find an adaptive law that will lead to a negative semi-definite $\dot V$. The efforts in \cite{Fom81,Ort96,Annaswamy1998,Loh_1999} developed a theory of adaptive control for nonlinearly parameterized systems. The resulting controller structure for the case when $f$ is concave/convex in $\theta$ is summarized below:
    \begin{align}
    u & = -f(X_p,\hat\theta)+\hat\alpha^T X_p+r -a^*s\left(\frac{e_c}{ \epsilon}\right)  \label{adaptive_control}\\
    s(y) &= \begin{cases}
     y^{(2\beta+1)} & \text{if} \ |y| < 1  \\
     sgn(y) & \text{otherwise}
     \end{cases} \label{s_def}\\
    a^* & = \lambda_{max} \min_{\omega \in \mathbb{R}^n} \max_{\theta_i \in \Theta_s} sgn(e_c)J\label{a_i_def}\\
    \omega^{*} & = \mbox{arg.}  \min_{\omega \in \mathbb{R}^n} \max_{\theta_i \in \Theta_s} sgn(e_c)J\label{omega_i_def}\\
    J & = \beta\left(f(X_p,\theta)-f(X_p,\hat\theta)+ \Tilde{\theta}^T\omega\right) \label{J_i_def}
    \end{align}
where $\beta$ is a known constant and $\Theta_s$ is a known compact set that the parameter $\theta$ belongs to. The update laws for the adjustable parameters in \eqref{adaptive_control} are chosen as:
    \begin{align}
    \dot {\hat{\alpha}} & = -\Gamma_{\alpha} e_{\epsilon}X_p \label{adcontrol2}\\
    \dot {\hat{\theta}} & = \Gamma_\theta e_{\epsilon}\omega^*\label{adcontrol3}
    \end{align}
where $e_\epsilon=e_c-\epsilon s(e_c/\epsilon)$, and $\Gamma_\alpha$ and $\Gamma_\theta$ are symmetric positive definite matrices. Closed-form expressions for $a^*$ and $\omega^*$ can be found when $\beta f(X_p,\theta)$ is convex for all $\theta\in\Theta_s$ or concave for all $\theta\in\Theta_s$.  While the solutions for $\omega^*$ coincide with the gradient of $f$ in some cases, they do not in other cases. Extensions to the case when $f(.,\theta)$ is a general function of $\theta$ can be found in \cite{Loh_1999}. Properties of persistent excitation that guarantees parameter learning have been addressed in \cite{Cao_2003}.

Three main points should be noted: Adaptive control approaches can be applied to problems where the underlying nonlinearities are convex (or concave). This makes the methodology applicable for nonlinearities that can be approximated by neural networks with convex activation functions such as ReLU \citep{Patkar2020}. The second point to note here is that even for these convex functions, new tools that are beyond the deployment of gradient methods such as min-max tools have to be introduced to lead to global solutions. The third point is that  powerful tools as Lyapunov functions based on Bregman divergence \citep{Bof21}, that allow a better accommodation of nonlinearly parametrized systems,  may still have a problem when  dealing with deep neural networks. This may be because of the latter introducing significant nonconvexities such as nonconvex dependence of the underlying loss function on weights of hidden layers.

\subsection{Stochastic and discrete-time Systems: Stability}\label{ss:stochasticstability}
The major milestone in adaptive control of stochastic and discrete-time systems is the proof of stability \citep{Goodwin_1981,Astrom_1995,Solo79,Landau82,Bit83,Kumar83,clarke1985generalized,Caines84}. We summarize this result by grouping various highlights in the literature under two headings: (1)	SA and RLS algorithms; (2)	Proof of stability of STR.

\subsubsection{Parameter estimation algorithms}
Several problems in system identification and adaptive control can be reduced to the identification of an unknown parameter vector $\theta^*$ in \eqref{eq:11} using input-output data stemming from regression vector $\phi_k$ and the output $y_k$. Two well known algorithms, recursive in nature, have been developed in the 70s and 80s and played a major part in adaptive control. These are described below:
\paragraph{Stochastic approximation based algorithm} 
Denoting  $\theta_k$  as the recursive estimate of $\theta^*$ as using an estimator, an estimated output for the system in \eqref{eq:11} can be derived as
\begin{equation}
  \hat{y}_k = \phi_{k-1}^\top\theta_{k-1}.
  \label{eq:estimate}
\end{equation}
The stochastic approximation (SA) algorithm takes the form \citep{Becker85}
\begin{align}
  \theta_{k} &= \theta_{k-1} -\frac{\gamma}{r_{k-1}}\phi_{k-1}(\hat y_k-y_k), \gamma>0 \label{eq:kumar3}\\  
  r_k &= r_{k-1}+\phi_k^T\phi_k;\;\;r_0=1\label{eq:kumar3a}
\end{align}
Several variations of the algorithm in \eqref{eq:kumar3}-\eqref{eq:kumar3a} have been proposed in the literature. Denoting $M_k=\gamma/r_k$, $M_k$ can be chosen to be a matrix rather than as in \eqref{eq:kumar3}. Instead of \eqref{eq:kumar3a}, a non-recursive choice of $r_k=1+\phi_k^T\phi_k$ can be introduced, which coincides with the projection algorithm in \cite{Goodwin_1981}. The following theorem summarizes the properties of the projection algorithm when there is no noise:
\begin{theorem} 
For the system in \eqref{eq:11} with $v_k\equiv 0$, it can be shown that
\begin{enumerate}
    \item $||\wh_k-\theta^*||\leq ||\wh_{0}-\theta^*||\qquad\forall \;k$
    \item $\lm{k}||\wh_k-\wh_{k-t}||=0$ for any finite $t$.
\end{enumerate}\end{theorem}

Similar results exist for the case when $v_k\neq 0$. The reader is referred to \cite{Becker85} and \cite{Goodwin_1981}(Chapter 8) for further details. 
\paragraph{RLS algorithm:} A simple variation of the adaptive gain $M_k$ in the SA algorithm leads to the well known recursive least squares (RLS) algorithm, summarized below:
\begin{align}
  \theta_{k} &= \theta_{k-1} -\Gamma_{k-1}\phi_{k-1}(\hat y_k-y_k),  \label{eq:4}\\  
  \Gamma_k &= \Gamma_{k-1}-\Gamma_{k-1}\frac{\phi_{k}\phi^T_{k}\Gamma_{k-1}}{1+\phi_{k}^T\Gamma_{k-1}\phi_k};\label{eq:4a}
\end{align}
It is easy to see that the RLS algorithm in \eqref{eq:4}-\eqref{eq:4a} is a matrix version of the SA algorithm above as well as the well known Robbins-Munro algorithm where $\Gamma_k$ is replaced by a scalar gain $\gamma_k$ satisfying additional conditions as in Eq. \eqref{robbinsmunro}.

Similar to the discussions of parameter estimation in the continuous-time case, here too, convergence of the parameter estimates to their true values is predicated on the persistent excitation of the regressor $\phi_k$. Formally this is stated as follows \cite{Anderson_1982}, and can be viewed as a discrete-time analog of Definition 1:
\begin{definition}\label{discrete:PE}
    A bounded function $\phi:N\rightarrow\mathbb{R}^N$ is persistently exciting (PE) if there exists $T\hspace{-.045cm}>\hspace{-.045cm}0$ and $\alpha\hspace{-.045cm}>\hspace{-.045cm}0$ such that
    \begin{equation*}
        \sum_t^{t+T}\phi_k\phi^T_kd\tau\geq\alpha I,\quad \forall t\geq 0.
    \end{equation*}
\end{definition}

\subsubsection{Adaptive control: Proof of stability}
As in the previous section, we state the main result for the noise-free case, and defer the reader to \cite{Kumar85,Goodwin_1981} for the noisy case. The starting point is the system in \eqref{kumar1}, with $w_k=0$. Defining polynomials $F(z)$ and $G(z)$ as 
\be 1= AF+z^dG\label{bezout}\ee it is easy to see that a control input $u_k$ chosen in the form of
\be u_k = - (G/FB) (y_k+r_k) \label{fixed-discrete}\ee
where $r_k=y^*_{k+d}$ for any bounded sequence $y^*_k$ ensures that $y_{k+d}=y^*_{k+d}.$ That is, the tracking problem is solved by choosing the control input in the form \eqref{fixed-discrete} when the parameters of $A$ and $B$ are known, provided the system is minimum-phase, that is, all roots of $B(z)$ are inside the unit circle. With this restriction, one can proceed to determine the stabilizing adaptive controller.

We reparameterize the polynomials $F$ and $G$  in the form
\bena 
F(z)B(z)&=&\beta_0[1+\beta(z)]; \beta(z)=\beta_1z+\ldots+\beta_{m+d-1}z^{m+d-1}\\
\frac{1}{\beta_0}G(z)&=&\alpha_0+\alpha_1z+\ldots +\alpha_{n-1}z^{n-1}\eena
and collect the coeffcients of $G$ and $FB$ as $\theta_c^*=[\alpha_0,\ldots,\alpha_{n-1},\beta_1,\ldots,\beta_{m+d-1},\frac{1}{\beta_0}]^\top$. 
We then write the system \eqref{eq:11} in a predictor form
\be u_k=\phi_{c,k}^T\theta_c^*\label{predict-form}\ee
where $\phi_{c,k} = [-y_{k}, \ldots, -y_{k-n+1}, -u_{k-1}, \ldots, -u_{k-m-d+1},y_{k+d}]^\top$.
This allows us to express the desired control input $u^*_k$ in \eqref{fixed-discrete}, when the parameters are known, as
\be u^*_k=\varphi_k^T\theta_{c}^*\label{fixed-discrete2}\ee
where $\varphi_{k} = [-y_{k}, \ldots, y_{k-n+1}, -u_{k-1}, \ldots, -u_{k-m-d+1},r_k]^\top$.
This leads to an adaptive controller
\be u_k=\varphi_k^T\theta_{c,k}\label{adaptdiscrete} \ee
with the parameter estimate $\wh_{c_k}$ adjusted using a variation of the SA algorithm
\begin{align}
  \theta_{c_k} &= \theta_{c_{k-1}} +\frac{\gamma}{c+\phi^T_{c,k-d}\phi_{c,k-d}} \phi_{c,k-d}(u_{k-d}-\phi^T_{c,k-d}\theta_{c_{k-1}}),  \label{eq:13}
\end{align}
where $0<\gamma< 2$ \citep{Goodwin_1981}. The following theorem summarizes the main stability result:
\begin{theorem}
Under the assumptions that (i) $d$, $n$, and $m$ are known, (ii) the zeros of $B(z)$ lie inside the unit circle, (iii) that there are no common factors between $A(z)$ and $B(z)$, and (iv)$\beta_0\neq 1$, the following hold:
\begin{enumerate}
    \item $\{y_k\}$ and $\{u_k\}$ are bounded sequences, \item $\lm{N}|y_k-y^*_k|=0$, and \item $\lm{N} \sum_{i=d}^N[y_i-y_i^*]^2<\infty$
\end{enumerate}
\end{theorem}
The result above establishes clearly that a minimum variance controller can be obtained when the underlying parameters of a system as in \eqref{eq:11} are not known. The adaptive algorithms can either be of SA-type as in \eqref{eq:13} or an RLS-type as in \eqref{eq:4}-\eqref{eq:4a}.

As in the continuous time case presented in Section \ref{solutions-stability}, an equivalent set of assumptions needs to be satisfied for the stability result here to hold. These correspond to the following: (i) the order $n$, and the delay $d$ have to be known. (ii) the sign of $\beta_0$ and a lower bound on the magnitude of $\beta_0$ need to be known. (iii) the zeros of $B(z)$ have to lie inside the unit circle. As outlined in Theorem 2, when these assumptions hold, a real-time adaptive control solution can be derived for the control input which ensures that for any initial conditions of the states and the parameter estimates, that the output error converges and is in $l_2$. Parameter learning follows as in the continuous-time case with persistent excitation of $\varphi_k$.

\subsubsection{Adaptive LQG control}
The problem of adaptive control when the underlying cost is quadratic both in the states and the inputs, as in \eqref{kumar11}, becomes much more difficult, and requires several more additional assumptions and results in weaker results. These are summarized below.

We return to the problem statement in Eqs. \eqref{kumar10} and \eqref{kumar11}. It is well known that for this linear-quadratic-guassian system the following control input is optimal:
\be u_k=K(A,B)x_k\label{kumar20}\ee
where 
$$K(A,B)=-[B^TPB+T]^{-1}B^TPA$$ and
$$P=A^TPA-A^TPB(B^TPB+T)^{-1}B^TPE+Q.$$
The results in \cite{Becker85,Campi96,Campi98} clearly show that the problem becomes significantly more complex when $A,B$ are unknown, and the control gain in \eqref{kumar20} has to be replaced with that which depends on parameter estimates of $A,B$. Suppose we define $(A_k^{LS},B_k^{LS})$ as the least squares estimate of $[A,B]$,i.e. 
\be
(A_k^{LS},B_k^{LS})\eq {\rm arg min}_{(A,B)\in\Theta}\sum_{s=1}^k||x_s-Ax_{s-1}-Bu_{s-1}||^2\label{kumar21}\ee
It is shown in \cite{Becker85} for ARMAX systems that the parameter estimates can converge to false values with positive probabilities; an example of the above statement for general Markov chains can be found in \cite{Borkar79}. 

A few interesting extensions have been reported in \cite{Campi96,Campi98} towards a suboptimal and stable solution under additional assumptions. This is accomplished by adding a bias term to the cost $J$ in \eqref{eq:11} so as to lead to estimates of the form  
\begin{align}
(A_k^{LS},B_k^{LS})&= {\rm arg min}_{(A,B)\in\Theta}\sum_{s=1}^k||x_s-Ax_{s-1}-Bu_{s-1}||^2\nonumber\\
&\qquad +\mu_tJ(A,B) \qquad {\rm if}\; k\;\hbox{is even}\label{kumar22}\\
&=(\widehat A_{k-1},\widehat B_{k-1})\qquad {\rm if }\; k \;{\rm is odd}\label{kumar23}\end{align}

\bluec{In addition to the above, the use of diminished persistent excitation with time was utilized to lead to adaptive optimal control in stochastic systems in \citep{GuoChen91,Guo95,Duncan99}.}
\subsection{Adaptive control of Continuous-time systems: Robustness}\label{sec:adapt-robust}
Suppose we start with an input-output model of an uncertain linear dynamic system \eqref{plant-linear}. The question that immediately arises is as to what is uncertain in \eqref{plant-linear}. The path that has been adopted in the field of adaptive control is to lump the uncertainty entirely into $\theta$ in \eqref{plant-linear}, the parameter of the dynamic system. The results outlined above, in Sections \ref{solutions-stability}-\ref{ss:stochasticstability}, proceeded  with such a problem statement as the starting point. The next step in the evolution of adaptive control expanded the scope of the problem from \eqref{plant-linear} to \eqref{plant-linear-robust}, where parametric uncertainties in $\theta$ were assumed to be accompanied by non-parametric uncertainties in the form of $d(t)$, $\vartheta(t)$, and $\Delta(s)$. The question that was addressed was how the solutions developed for \eqref{plant-linear-robust} can remain satisfactory even with these non-parametric perturbations. Two broad classes of solutions were proposed in the literature, one that sought to modify the adaptive controller in \eqref{ad1}-\eqref{ad2}, by changing the adaptive law in \eqref{ad2} to a form
\bea \dot\theta_c &=& C_2(\theta_c,\phi,t) -h(\theta_c,\phi)\label{robustadaptive1} \eea
where the correction term $h(\cdot)$ is designed to produce robustness. The second type of results retained the adaptive control structure as in \eqref{ad1}-\eqref{ad2} but invoked conditions of persistent excitation on the exogenous signal. These are summarized in the following sections.

\subsubsection{Modifications in the adaptive law}
For ease of exposition, we restrict our discussion to linear systems with single-input, whose states are accessible. The reader is referred to textbooks such as \citep{Narendra2005,Ioannou1996} for further details.
Consider a plant-model of the form 
\be \dot x = A x + bu + v\label{robust-adaptive1}\ee where $A$ is an unknown matrix, $b$ is a known vecor, $u$ is the control input and is a scalar, $x$ is the state and is accessible for measurement, $v$ is a disturbance that is unknown, time-varying, and bounded. The use of matching conditions \eqref{e:match} suggests that a reference model of the form
\be \dot x_m=A_m x_m+ br\label{e:ref-model}\ee
where $A_m=A+bk^{*T}$ is known and Hurwitz generates a class of command signals $x_m(t)$ that the plant state can be guaranteed to track, by choosing a control input $u=k^{*T}x+r$. As $A$ is unknown, an adaptive control input and adaptive law of the form
\be 
u=k^Tx+r,\qquad \dot k = - \gamma (e^TPb)x \label{adapt-state}\ee
where $\gamma>0$ is a positive constant, $e=x-x_m$, and $P$ solves the Lyapunov equation \eqref{Lyap} guarantees that 
\be V=e^TPe+(k-k^*)^T(k-k^*)\label{lyap-state2}\ee
is a Lyapunov function with $\dot V=-e^TQe$ and that $\lm{t}e(t)=0$, provided $v(t)\equiv 0$. When $v(t)\neq 0$, the same stabilizing control input in \eqref{adapt-state} contributes to a parameter drift in $k$ to $-\infty$ \citep{Roh85,Nar86}. This is because of a windup effect and the fact that the adaptive controller is a nonlinear integral controller; in the presence of a disturbance, it can cause the parameter to wind-up to infinity. 

The solutions suggested in the literature introduce anti-windup actions in the form of a correction to the adaptive law in \ref{adapt-state} as
\be 
\dot k = - \gamma (e^TPb)x -h(e,x,k)\label{adapt-state2}\ee
which causes the time-derivative to take the form
\be \dot V = -e^TQe +2 e^TPv - 2k^T h(k,e,k) \label{lyap-dot-robust}\ee
The approaches in the literature pertain to different choices of $h(e,x,k)$ such that $\dot V<0$ outside a compact set in the $(e,k)$ space. An equivalent approach is to modify the underlying loss function such as the one in \eqref{e:error2_Loss} with a regularization term that involves the ${\cal L}_2$ norm of $k$. 
It should be noted that existing literature \citep{Narendra2005,Ioannou1996} includes results for the case when partial set of states are available for measurement, when there are multiple inputs or when the underlying system is in discrete-time \citep{disIoa86,Tao95,Wen92,Clu88}.

The discussions above were focused on the perturbed model in \eqref{plant-linear-robust} where there is either a disturbance $d(t)$ or the parameter $\theta$ is a function of time. Robustness to unmodeled dynamics such as $\Delta(s)$ in \eqref{plant-linear-robust} is a considerably more challenging problem compared to robustness to either bounded disturbances or time-varying parameters. Of equal difficulty is robustness to time delays, which are ubiquitous in large-scale and networked systems. Several results have been proposed in the context of robustness to unmodeled dynamics (for example, \cite{Narendra2005,Ioannou1996,Naik92,Pomet92} and time-delays (for example, \cite{Ort88,Nic03,Yil10,Bresch2009}), and more recently in \cite{Hus17,Hus17thesis,Dog16,Hus13}) many of which employ a projection operator \citep{Lav2011}. These establish that adaptive systems can be designed to be robust with respect to unmodeled dynamics by having the parameters adapt inside a bounded set and guarantee bounded solutions.

\subsubsection{Use of Persistent Excitation}
An alternate approach to establish robustness, i.e. bounded solutions in the presence of the disturbance $v$ in the adaptive system defined by \eqref{robust-adaptive1}-\eqref{adapt-state} is to invoke conditions of persistent excitation of $x_m$ in \eqref{e:ref-model}. Two classes of results have been reported in the literature, in \cite{Nar86} and in \cite{And86} using such an approach. The results in \citep{Nar86} are briefly summarized below for the adaptive system in \eqref{robust-adaptive1}-\eqref{adapt-state}. Let us assume that the exogenous input $r$ is such that $x_m(t)$ is persistently exciting in $\RR^n$ with the level of persistent excitation $\epsilon_0$, which is defined as \ben |\frac{1}{T_0}\int_{t_2}^{t_2+\delta_0}x_m^T(\tau)w d\tau|\geq \epsilon_0 \;\forall\; t\geq t_0 \een 
where $[t_2,t_2+\delta_0]\subset [t,t+T_0]$ and $w$ is a unit vector in $\RR^n$. Then the adaptive system will have globally bounded solutions if
\ben \epsilon_0 > k_0 v_{\rm max} \een where
$|v(t)|\leq v_{\rm max}$ and $k_0=2\lambda_{P-max}/\lambda{Q-min}$. That is, if the level of persistent excitation is large compared to the size of the disturbance, then boundedness follows. It is also shown in \cite{Nar86} that the converse is true - for a class of adaptive systems, for a class of disturbances, it can be shown that there exists a signal $x_m(t)$ for which solutions of the adaptive system will be guaranteed to exhibit instability in the form of $\lm{t}k(t)=-\infty$. A similar phenomenon was shown in \cite{Roh85} to hold in numerical simulations.
The results of \cite{And86} established a similar result for the harder problem when $v(t)$ is not necessarily bounded, but state-dependent, which occurs when it is due to unmodeled dynamics excited in closed-loop. The authors therein showed that when the underlying regressors are persistently exciting, the properties of the adaptive system can be locally approximated by an averaged system that has exponential stability properties and therefore shown to be robust. 

\subsection{Bursting Phenomenon and Imperfect Learning}\label{ss:bursting}
The results in the above section clearly indicate the close relationship between the trajectories that the parameter estimates take,  persistent excitation, and disturbances. We point out another interesting property that has been observed in the context of adaptive systems and learning, which is the bursting phenomenon \citep{And86}. The milestones above indicate three distinct facts: (1) Persistent excitation of the underlying regressor leads to parameter convergence \citep{Morgan_1977,Anderson_1982}; (2) Persistent excitation at a sufficient level relative to the disturbance ensures robustness \citep{Nar86,And86}; (3) When the excitation level is not sufficient or if there is simply no persistent excitation, then parameters will not converge to the true values \citep{Kumar83}, i.e. leads to imperfect learning. A fourth fact that rounds off this topic is this: (4) When there is no persistent excitation, and when there are disturbances present, the closed-loop system can produce large bursts of tracking error \citep{Mor77,For81,And85}. That is, imperfect learning exhibits a clearly non-robust property that leads to a significant departure from a tracking or a regulation goal: exhibit an undesirable behavior over short periods during when the tracking error becomes significantly large. 

A specific example that illustrates this behavior is the following \citep{And85}:
Consider a first-order plant with two unknown parameters $a$ and $b$ of the form
\be y_{k+1}=a y_k+b u_k \label{And1}\ee
whose adaptive control solution is given by \cite{Goodwin_1981}
\be u_k=\frac{1}{\widehat b_k}\left[-\widehat a_k y_k+ y^*_{k+1}\right]\label{And2}\ee
The results of \cite{Goodwin_1980} in \eqref{adaptdiscrete} and \eqref{eq:13} reparameterize \eqref{And2} as 
\be u_k=-\theta_{c_1,k} y_k + \theta_{c_2,k}y^*_{k+1}\label{And3}\ee
and propose a parameter adjustment rule as in \eqref{eq:13} where $\phi_{c,k}=[-y_k,y_{k+1}]^\top$ and $d=1$. Clearly, the results in the literature  guarantee that the adaptive controller defined in \eqref{And3},\eqref{eq:13} guarantee that (i) $\theta_{c_i,k}$ and $y_k$ are bounded \citep{Goodwin_1981} (ii) $\theta_{c_i,k}$ converge to constants $\theta_{c_i}^0$, which may not coincide with the true values \citep{Becker85}, and that (iii) $y_k$ approaches $y^*_k$ as $k\tends\infty$ \citep{Goodwin_1981}. In addition, when $\phi_{c,k}$ is persistently exciting, i.e., satisfies Definition 2, 
we also have that the estimates $\theta_{c_i,k}$ approach the true values $\theta_{c_i}^*$. When such a persistent excitation is not present and when perturbations are present, bursting can occur, which can be explained as follows:

Suppose we consider a simple regulation problem with $y^*_k\equiv 1$. The control input in \eqref{And3} leads to a closed-loop system of the form
\be y_{k+1} = g(\theta_{c_1,k})y_k+h(\theta_{c_2,k})\label{And4}\ee where
\be g(\theta_{c_1,k})=\left(a-b\theta_{c_1,k}\right),\qquad h(\theta_{c_2,k})=b\theta_{c_2,k}\label{And5}\ee
This implies that the closed-loop system is (a) unstable if $|g(\theta_{c_1,k})|> 1$, and (b) stable if $|g(\theta_{c_1,k})|< 1$. 
The most troublesome scenario occurs if $\theta_{c_1,k}=\theta_{c_1}^b$ where $g(\theta_{c_1}^b)=- 1$. Such a case will cause bursting.
When disturbances are present, the discussions in Section \ref{sec:adapt-robust} showed that parameters can drift. It is therefore possible that  parameters $\theta_{c_i,k}$ become arbitrarily close to $\theta_{c_i}^b$ for some $k=k_0$; at $k_0^+$ a disturbance pulse is introduced, causes the parameters to drift with $\theta_{c_1,k}$ approaching $\theta_{c_1}^b$, which in turn causes $y_k$ to oscillate, which then causes $\theta_{c_i,k}$ to readjust, once again approach another set of constant values $\theta_{c_i}^{0\prime}$. Such a phenomenon has been shown to occur in \cite{And85} and in continuous-time systems \citep{Nar87}. It should be noted that this occurs with imperfect learning, that is, when the underlying regressors are not persistently exciting. Such a phenomenon is not peculiar to the specific systems in question, but for any arbitrary dynamic systems where simultaneous identification and control are attempted. 

\subsection{Adaptive Control in the Presence of Input and State Constraints}\label{ss:constraints}
The adaptive controllers outlined in Sections  \ref{solutions-stability}-\ref{ss:stochasticstability} were focused on ensuring that the closed-loop system has bounded solutions and that the output error was minimized. No restrictions were imposed on the requisite control input. A wider problem statement with the goal of stable adaptive control in the presence of magnitude and rate constraints on inputs and states  was addressed in a number of publications including \citep{karason1994,Lav04,Hov08,Leon09,Gaudio_2019b} and is briefly summarized below.

Suppose that the output of the adaptive controller is denoted as $u(t)$, and the actual input into the plant is denoted as $u_p(t)$. Suppose that an elliptical saturation function is defined as 
$E_s$ denotes an elliptical saturation function of a vector $v(t)$ defined as \citep{Gaudio_2019b}
\begin{equation}\label{e:Es}
E_s(v(t),v_{max}) = 
\begin{cases}
v(t), &\quad||v(t)||\leq g(v(t))\\
\bar{v}(t), &\quad||v(t)||>g(v(t)) \ 
\end{cases}
\end{equation}
where the function $g(v(t))$ is expressed as
\begin{equation}
g(v(t))=\left(\sum_{i=1}^{m}\left[\frac{\hat{e}_i}{(v_{max})_i}\right]^2\right)^{-1/2},
\end{equation}
where $\hat{e}=\frac{v}{||v||}$ and $\bar{v}=\hat{e}g(v)$. The plant input $u_p$ is then generated using $E_s$ as
\bea 
u_r(t)&=& \frac{1}{\tau}(E_s(u(t),u_{max})-u_p(t))\label{e:sat1}\\
\dot{u}_p(t)&=&E_s(u_r(t),u_{r,max}) \label{e:sat2}
\eea
Such an input $u_p(t)$ is guaranteed to meet the magnitude limit, with $|u_p(t)|\leq u_{max}$, and coincides with $u(t)$ when the magnitude of $u(t)$ is small. In terms of rate, first we note that $|u_r(t)|\leq u_{r,max}$. In addition, since the variable $u_r(t)\approx \dot u$, it  follows that  $u_p(t)$ is rate limited with a bound of $u_{r,max}$. The reader is referred to the Appendix and \citep{karason1994,Gaudio_2019b} for details.

The saturation functions in magnitude \eqref{e:sat1} and rate \eqref{e:sat2}  introduce two nonlinearities, which poses a problem in the overall analysis. The main idea articulated in \citep{karason1994,Hov08,Gaudio_2019b} that overcomes this problem is to represent them as additive known disturbances. In particular, defining two known disturbance terms $\Delta u_m$ and $\Delta u_r$ as
\begin{align}
\label{e:Control_Deficiency}
\begin{split}
\Delta u_m(t)&=E_s(u(t),u_{max})-u(t),\\
\Delta u_r(t)&=E_s(u_r(t),u_{r,max})-u_r(t),
\end{split}
\end{align}
it is easy to see that if $u$ does not reach its magnitude saturation limit $u_{max}$, then $\Delta u_m(t)\equiv 0$. Similarly, if the input rate $u_r$ does not reach its rate saturation limit $u_{r, max}$, then $\Delta u_r(t)\equiv 0$, that is, these known disturbance terms become non-zero only if the magnitude or rate limits are exceeded.

Using (\ref{e:sat1}), \eqref{e:sat2}, and (\ref{e:Control_Deficiency}) we obtain a compact relation between the plant input and the controller output of the form
\begin{equation}
\label{e:filter}
\dot{u}_p(t)=-\frac{1}{\tau}u_p(t)+\frac{1}{\tau}u(t)+\frac{1}{\tau}\Delta u(t),
\end{equation}
where $\Delta u(t)=(\Delta u_m(t)+\tau\Delta u_r(t))$ represents the combined effects of magnitude and rate saturation.
That is, both magnitude and rate limits can be accommodated in the form of an additive disturbance $\Delta u(t)$ and a filter $1/(\tau s+1)$. 
This in turn implies that an underlying plant model of the form 
\be y_p=W_p(s)u_p \ee 
where the input $u_p$ is subject to magnitude limits and rate limits can be rewritten as
\be y_p=\frac{W_p(s)}{\tau s+1} (u+\Delta u)\label{e:sat3}\ee 
The effect of the magnitude limit is in the form of a disturbance while that of rate limit is in a combined form of both a filter and a disturbance. The adaptive control solutions propose in the literature address the problem of determining the control input \eqref{e:sat3} when the parameters of $W_p(s)$ are unknown. We briefly describe the solution for the simple case when states are accessible and only magnitude limits are imposed.

We start with the problem statement in Section \ref{ss:MRAC} 
\be \dot x = A_p x + b_pu_p \label{e:state2}\ee
where $u_p$ is the plant input and is required to meet a hard magnitude constraint $u_{{\rm max}}$. Using the procedure described above, it is easy to show that
$$u_p=u+\Delta u$$
where $u$ is the output of an adaptive controller derived as in \eqref{adaptive1} and 
\be
\Delta u = \left\{\begin{array}{cc} 0 & {\rm if}|u(t)|\leq u_{{\rm max}}\\
u-u_{{\rm max}}& {\rm if}|u(t)|\geq u_{{\rm max}}\end{array}\right.\label{magsat}\ee
This in turn leads to an error model of the form 
\be \dot e = A_m e + b_p(\wt^T\omega+\Delta u)\label{magsat2}\ee
where $\wt=[(theta-\theta^*)^T, (k-k^*)]^T$, and $\omega=[x_p^T, r]^T$.
As $\Delta u$ is a known disturbance, an augmented error $e_u$ is generated as $e_u=e+e_a$, where 
\be \dot e_a = A_m e_a + b_mk_s(t)\Delta u\label{magsat2b}\ee
which includes an additional adjustable parameter $k_s(t)$. This in turn allows an error model to be derived in a standard form as
\be \dot e_u = A_m e_u + b_p\wt^T\omega+b_m (k_s(t)-k^*)\Delta u\label{magsat3}\ee
A Lyapunov function similar to \eqref{lyap-state} can be found that guarantees the boundedness of $e_u$ and the adaptive parameters $\theta$, $k$, and $k_s$. An additional and significant hurdle now needs to be overcome to show boundedness of the plant states, as it can no longer be concluded that the original state error $e$ is bounded, as $e_u$ is a sum of two signals, both produced by the control input in closed-loop. Properties of linear systems with bounded inputs are employed in \citep{karason1994,Hov08,Gaudio_2019b} in order to show global boundedness for all open-loop stable plants and boundedness in a domain of attraction otherwise. Similar results have been derived for discrete-time plants as well in \citep{Kar95,Zha87,Chao01}.

\subsubsection{State Constraints and Barrier Functions}
Novel extensions to nonlinear systems with state constraints have been addressed in the literature through the use of Barrier Lyapunov functions \citep{Tee2009,ren2010,ames2014}. The main idea here is to construct Lyapunov functions that become large when the error variables approach certain limits. For example, rather than choose a quadratic term in $x$, a log function of the form $log (x^2_0/(x_0^2-x^2))$ is utilized to make sure that the state variable $x$ does not exceed its limit $x_0$.

\subsection{Assumptions and Challenges}\label{ss:assn}
The solutions outlined in Sections IV-A through IV-H correspond to decision-making in a dynamic system in real-time by a controller. As shown in the schematic in Fig. \ref{fig:proj}
\begin{figure}[!htb]
  \centering
  \includegraphics[width=0.5\textwidth]{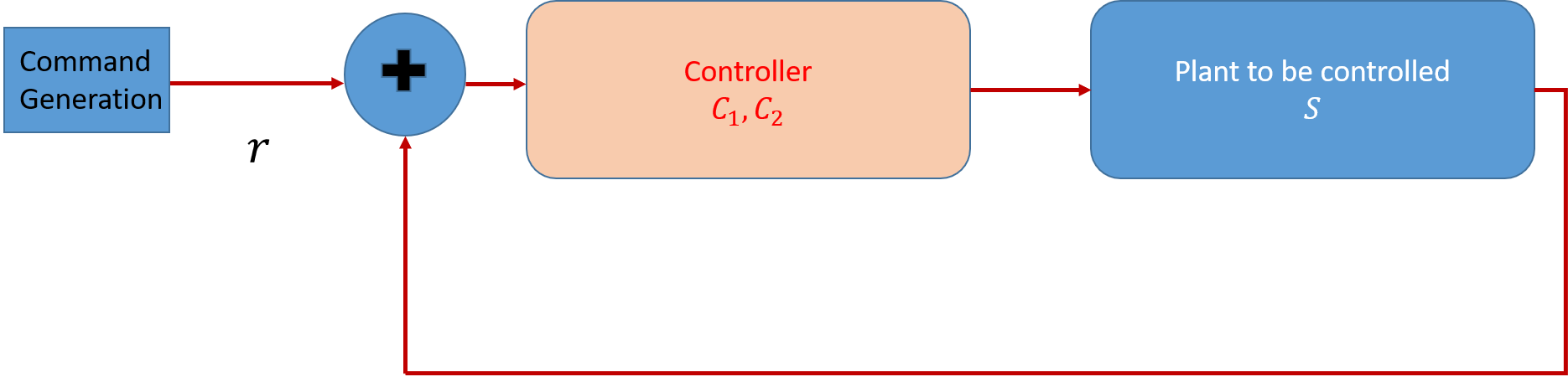} 
  \caption{Closed-loop control of a nonlinear dynamic system using adaptive control}
  \label{fig:proj}
\end{figure}the dynamic system $S$ is described  as in \eqref{plant-linear}, \eqref{kumar1}, \eqref{robust-adaptive1} for the linear case, or \eqref{plant1}, or \eqref{plant} in the nonlinear case. The adaptive controllers $C_1,C_2$ have a general form as in \eqref{ad1}-\eqref{ad2}, which in linear systems are of the form \eqref{e:Gradient_Flow_UN} in simple cases, \eqref{adaptive1}-\eqref{adaptive3} when states are accessible, \eqref{e:adcon} and \eqref{ad20}-\eqref{ad21} or \eqref{adlaw-general} for adaptive output feedback. For nonlinear problems, a few examples were outlined in \eqref{adaptive_control}-\eqref{adcontrol3}. In most of the cases, it should be noted that the solutions provided are global, with the adaptive system starting from arbitrary initial conditions, and are applicable in real-time. No training, exploration or simulation experiments are required. These are unique features and advantages of the adaptive control method. The guarantees that the adaptive control solutions provide are predicated on assumptions that some prior information is available about the plant. Examples are information about its  order,  net-order, and the sign of the high frequency gain for the linear case. In a nonlinear plant, the assumptions pertain to a certain type of interconnection such as strict-feedback form \citep{Krstic_1995}, triangular structures \citep{Seto_1994}, or feedback linearization \citep{Slotine1991}. It should again be mentioned that in almost all of these cases, the underlying solutions are provably correct, with firm analytical guarantees and precise descriptions of the nature of the solutions.

Several efforts have been consistently and continuously applied over the years to relax these assumptions. An approach credited to Nussbaum \citep{Nuss83} relaxes the requirement that the sign of the high frequency gain be known. The approach outlined in section \ref{s:general} under high-order tuners only requires the net-order, but not the order, to be known. Several extensions to  nonminimum phase systems have been reported over the years \citep{Elliott85}. In some cases, these extensions come with other drawbacks such as lack of robustness due to an intrinsic high-gain nature, requirements of persistent excitation, or significant increase in the computational burden.

The importance of the interconnection between adaptation and learning cannot be over emphasized. The stability results outlined in Section \ref{solutions-stability} focused on just that, stability. No requirements on learning the unknown parameter were involved. Additional conditions of persistent excitation, if imposed, led to learning. In closed-loop, as one cannot guarantee that such PE conditions can be satisfied, one cannot guarantee learning. Rather, with imperfect learning, the adaptive control solutions guaranteed closed-loop boundedness and asymptotic guarantee of performance. If additional conditions are met, then learning follows. Optimality of the adaptive controller, when learning is complete, can subsequently be ensured. As was evident from the discussions in Section \ref{sec:adapt-robust}, robustness of these performance goals under various conditions and perturbations have also been guaranteed.

As we proceed to expand the scope of the class of dynamic systems under consideration, we articulate one of the main challenges that is encountered due to the interconnection and compelling tradeoffs between adaptation and learning. Suppose we address the closed-loop control in Figure \ref{fig:proj}, where $S$ is given by \eqref{plant1}. 

Suppose that the control input $C_1$ is generated using a neural network with its input $\omega$ as follows:
\bea
y_i &=& \sum_{j=1}^{N_i}W^{iT}_{2j}\phi\left(W^i_{1j}y_{i-1}+b_{j}^i\right),\;\; i=2,\ldots L 
\label{nnet}\eea
where $L$ denotes the number of layers, $W^i_{kj}$, $k=1,2$, $b^i_j$, $j=1,\ldots N_i$ denote the weights in the $i$th layer, $y_1=\omega$, and $y_N=u$ are respectively the input and output of the neural controller.  $\phi$ denotes the activation function.  Suppose that the neural network is trained extensively for a given set of parameters in the dynamic system so that the closed-loop system generates a performance that is satisfactory. It should be noted that the parameters of the neural network, $W^i_{kj}$ converge (if they do; no analytical guarantees exist as of today that they converge) to some value $W^{i0}_{kj}$ such that the mapping between $\omega$ and $u$ approximate the desired nonlinear mapping. 
At this point, suppose that the plant parameter $\theta$ in $S $ changes in an unbeknownst manner that cannot be anticipated beforehand and  accompanied by non-parametric changes in disturbances or unmodeled dynamics. The closed-loop system in such a case is highly prone to the bursting phenomenon described in Section IV-I, in the multi-dimensional space made up of weights, states, and inputs of the overall closed-loop system, as the system may not have been trained satisfactorily at these changed conditions. It is the analysis of this resulting closed-loop system together with guarantees of boundedness, convergence, and optimality that is needed.  This is an open problem that needs to be addressed.

\subsection{Loci of Adaptive Control}

\bluec{As mentioned in the introduction, the focus of this article is on those aspects of adaptive control that has an identifiable learning component with tractable problem formulations and solutions. It should be noted that huge swaths of efforts have been expended in several other branches of adaptive control  over the past five decades  with enormous success. We mention but a few of those classes that capture the loci of adaptive control. }

\bluec{Extensions to adaptive control of infinite dimensional systems can be found in \citep{Smys2010} and a special class of problems which corresponds to systems with delays \citep{Ort88,Nic03,Yil10c,Bresch2009,bresch2014}. They have found applications in traffic control \citep{Burk21}, power-train control \citep{Yil10}, rocket pressure control \citep{Yildiz18}, and drilling \citep{krs2013}, to name a few.  In most of these cases, these approaches consist of adaptive controller designs with an in-built parameter estimate, with requisite complexities in both the control and adaptive law as well as in the machineries employed. Additional tools from infinite dimensional systems and Lyapunov functions need to be utilized to derive stable solutions. The goals of these designs are to primarily accomplish the control goal; learning, i.e. convergence of the underlying estimates to the true value are very hard to establish. }

\bluec{Interesting extensions have been reported in \citep{Guo97,XieGuo00,HuangGuo} for necessary and sufficient conditions for control of  classes of systems under uncertainties, in an attempt to examine fundamental limitations of the feedback mechanism. A computationally tractable solution to the adaptive stabilization problem addressed in \cite{HuangGuo} is proposed in \citep{Sokolov16a,Sokolov16b} based on  set estimation with Yakubovich's method of recursive goal inequalities. This and related directions are surveyed in \citep{Guo20}.}

\bluec{Several other branches of adaptive control  have been investigated over the years. The first include decentralized and distributed adaptive control based on notions of cooperation and consensus, and adaptive control for synchronization of complex networks \citep{ioannou86,DeLellis09,Olgren04,Zhou06,Cao08,Hou09,Das10}. Another area is adaptive control in the presence of commonly present algebraic nonlinearities such as hysteresis \citep{tao1995hyster} and deadzones \citep{tao1994dead}, which are useful in all applications where actuator nonlinearities have to be contended with \citep{tao2004actuator}. Along with actuator nonlinearities, actuator redundancy has been addressed in \citep{Yildiz20} via adaptive control allocation methods. The  use of a filter and high-gain in closed-loop \citep{Hov2010} has been explored as well with significant successes reported in applications.  Yet another related topic that intersects with adaptive control and machine learning based optimization is extremum seeking. Here, the goal is to adjust a parameter, but not with the purpose of learning the parameter or an underlying function but to rather maximize a function \citep{krstic2000extremum,ariyur2003extremum}. We expect interesting discoveries related to the intersections between these topics to unfold over the coming years. }




\bluec{The topics covered in this paper are by no means an exhaustive account of all control methods are adopted for dynamic systems with uncertain parameters. The most notable methods that we have not covered in this paper include adaptive sliding mode control (see for example, \cite{Bar95,Hua08,Lee09}), iterative learning control (see for example, \cite{Bris06},\cite{Moore12}), and linear-parameter-varying systems (see for example, \cite{Moham12,Hoff14}). The reader is referred to the cited papers for a deeper dive into these methods. }




\section{Applications}
Progress in theory has been accompanied throughout the past five decades with explorations of applications of adaptive control in various sectors. This is evidenced by edited books \citep{Nar80}, surveys \citep{Astrom83,Astrom96}, chapters in textbooks \citep{Landau11,Narendra2005,Uls89}, or entire textbooks \citep{Lavretsky2013}. Applications span process control \citep{Dum02,Dum03,Dum90,Dum95,Dum03}, automotive systems \citep{Yil10,Yil10b}, positioning systems \citep{Uls89,Smi95}, propulsion systems \citep{Eve03b,Ril04}, and a huge effort in flight control (see for example, \citep{Thom70},\citep{Tay64,Dyd10,Cal98,Jen00,Bos04},\citep{Gre11},\citep{Hol2011},\citep{Dydek12,Dyd13,Dydek13}).  

Since the early 1990’s, the US Air Force, US Navy, and NASA
working with industry and academia have made significant
progress towards maturing adaptive control theory for
aerospace applications \citep{Gre11}. Several adaptive control architectures
have been implemented in unmanned flight platforms
\citep{Sharma06}. A specific observer-based adaptive control with Loop
Transfer Recovery (OBLTR) has been developed in \cite{Lavretsky2013}. As
pointed out in \citep{Wise18}, a technology transition of conventional
MRAC applications and adaptive OBLTR based architectures,
has been continuously ongoing (see Figure 1, Figure 4 in
\citep{Wise18}) which includes aerial platforms such as JDAM, X-36,
and several others.

The reader is referred to the surveys and textbooks mentioned above for several more applications in addition to all of the above, for autopilots for ships, and coworkers, cement mills, chemical reactors, diesel engines, glass furnaces, heating and ventilation, motor drives, paper machines, optical telescopes and titanium oxide kilns, and more. The reader is referred to \cite{IoCT11} for additional success stories.

Several industrial products exist that implement MRAC and STR described above. NOVATUNE and NOVAMAX produced by ASEA AB were early examples mentioned in the 80s in \cite{Astrom83}. The reference \cite{Dum02} lists BrainWave, an adaptive MPC, and MicroController 2000X, both implemented in several process control problems. We also refer the reader to the proceedings of several workshops in conferences such as the ACC, CCTA, and CDC, that have presented recent applications of adaptive control to aerospace problems (see for example \citep{Hea21}).

In much of these applications, the need for adaptive control stems from a scenario where a control problem arises and a satisfactory solution requires a retuning of the control parameters due to aging, drift, or other untoward changes in the plant being controlled. The existing baseline controller becomes, as a result, incorrect, and needs to be retuned. In several of these applications, such as in autonomous vehicles either in air, ground, or water, may require this self-tuning or adaptation, to occur on the fly, in real time. The flight platforms considered in \citep{Sharma06} and other flight platforms listed above fall under this category. Often the existing baseline controllers become inadequate under these anomalies, may be destabilizing, and need to be retuned. And under these circumstances, adaptive control enables a procedure by which real-time adjustment of controllers is possible. That adaptive controllers are finding a pathway for technology transition, systematic validation, and field implementation is clear from the above discussions. The relatively slower pace of implementation of adaptive technologies, compared to say MPC, may be a combination of the need for a truly real-time tuning in a given application and the requisite bandwidth and complexity for implementation. As applications become more complex and as computing and communication technologies become more advanced, both of this impediments may very likely diminish and disappear.

\section{Summary and Concluding Remarks}
In this paper, we have sketched a historical perspective of the field of adaptive control over the past seven decades. Given the recent upsurge of interest in learning, both offline and online, in the Machine Learning and control communities,  such a perspective is timely and warranted. The scope of this article is large - we have attempted to cover highlights of the field which span  70 years, chronicled in $\sim$15 textbooks, $\sim$20 edited books, hundreds of surveys, and thousands of research publications in journals and conferences in 30 pages, which is a formidable task. We have therefore showcased just the highlights of this field, and emphasized key lessons learned, problems that have already been solved, important takeaway messages, and cautionary remarks. \bluec{While our attempts at chronology span the footprint of this topic  from the 1950s to the present, it should be acknowledged that there is a large vigorous set of activities in this area over the last five to ten years, especially at the intersection of parameter learning, reinforcement learning, neural networks, and adaptive control that we have not addressed in this paper. We refer the reader to recent plenary talks, papers in recent control and machine learning conferences, and special issues in related journals for the exposition of the latest advances.}

Over the last seventy years, the field of adaptive control has witnessed advances in both deterministic and continuous-time systems and stochastic discrete-time systems. We have attempted to cover both domains in this article. Key advances in different parts of the globe have all been attempted to be covered. We have not offered a deep technical discussion of theorems, but rather the idea behind key results and their implications. No proofs have been provided either. The reader is referred to the list of copious references at the end of the paper for in-depth technical expositions of all problems and solutions outlined here. We presented a chronological taxonomy of the advances in the field in Section II, a cross-section of problem statements in Section III, and highlights of key solutions in Section IV. Major applications of adaptive control are addressed in Section V.

The primary focus of the adaptive controllers has been to ensure that (a) the closed-loop system have bounded solutions, and (b) asymptotic properties of the outputs (and in some cases inputs) are achieved. The results in \ref{ss:constraints} have extended this focus and have made inroads in making sure that the  requisite constraints of magnitude and rate for the control input and state constraints are met as well. It should be noted that in all cases, the performance goals have been limited to the system behavior in real-time, at time $t$, and not for all future instants. As the premise in all these problems is that parametric uncertainties can be introduced at any time, optimization of a cost function over all time, with a cold-start of the controller that simultaneously estimates, adapts, and optimizes, is difficult if not impossible. Some of the recent results that propose clever combinations of both adaptive control and machine learning concepts may overcome this formidable challenge.

\section*{Acknowledgements}
We would like to thank P.R. Kumar for several useful discussions and directing us to the highlights of stochastic adaptive control. We would like thank the reviewers and Miroslav Krstic for their valuable comments which helped place our overall message in the broader and rich landscape of adaptive control theory. We would like to gratefully acknowledge Yingnan Cui and Boris Andrievsky
for helping us build the database with almost 300 references on adaptive control and learning.

\section*{Appendix}\label{sec:tools}
\subsection*{Stability framework}
\label{ss:Stability}
The first and foremost challenge introduced by adaptive control is a nonlinearity. As the controller is proposed as a real-time control solution, the nonlinearity is introduced due to the simultaneous estimation and control. That is, the control input is a function of the parameter estimate which in turn depends on the control input as well as several other system variables. As a result, the closed-loop system becomes nonlinear, with its solutions corresponding to the true responses of the plant being controlled. As a result, the well behavedness of the overall adaptive system, i.e. its stability is the first property that needs to be assured. The typical tool employed here is due to Lyapunov and is summarized below \cite{Narendra2005}. The dynamic system of interest is of the form 
\be \dot x = f(x,t), \qquad f(0,t)\equiv 0 \label{lyap1}\ee It is assumed that $f:\RR\tends\RR^n$ is such that a solution $x(t;x_0,t_0)$ exists for all $t\geq t_0$.
\begin{theorem} 
The equilibrium state $x=0$ of \eqref{lyap1} is uniformly asymptotically stable in the large if a scalar function $V(x,t)$ with continuous first partial derivatives w.r.t $x$ and $t$ exists such that $V(0,t)=0$ and if the following conditions are satisfied:
\begin{enumerate}
    \item[(i)] $V(x,t)$ is positive-definite,\item[(ii)] $V(x,t)$ is decrescent, \item[(iii)] $\dot V(x,t)$ is negative-definite, and \item[(iv)] $V(x,t)$ is radially unbounded.
\end{enumerate}\end{theorem}
$V(x,t)$ that satisfies these conditions is referred to as a Lyapunov function. If instead of (iii), $\dot V(x,t)$ is only negative semi-definite, only uniform stability can be ensured; if instead of (iii), a stronger condition $\dot V(x,t)< -\alpha(||x||)<- \beta V(x,t)$, then exponential stability of the equilibrium can be ensured. We refer the reader to \cite{Narendra2005} for all further technical details.

A typical approach in adaptive control is to express the underlying system in the form of \eqref{lyap1} with the state $x$ corresponding to errors in the system that are either to be driven to zero or required to be bounded. These error can be broadly grouped into two categories, tracking error and parameter error. Often adaptive systems consider quadratic Lyapunov function and only lead to a negative semi-definite $\dot V(x,t)$.
\subsection*{Rational SPR functions and the KYL}
The definition of SPR and one of the simplest versions of the KYL \cite{Anderson_1982} is given below.

Definition: An $n\times n$ matrix $Z(s)$, whose elements are rational transfer function, is SPR if for some $\epsilon>0$ and all $Re[s-\epsilon]>0$,\begin{enumerate} \item all elements of $Z(s-\epsilon)$ are analytic  \item $Z^*(s-\epsilon)=Z(s^*-\epsilon)$, and \item $Z^T(s^*-\epsilon)+Z(s-\epsilon)$ is positive semi-definite.
\end{enumerate}

\subsubsection*{The Kalman Yakubovich Lemma}
 Let $Z(s)$  be a matrix of rational functions with $Z(\infty)=0$, a minimal realization $\{A,B,C\}$, and with all its poles only in $Re[s]<-\mu$. Then $Z(s)$ is SPR if and only if there exist symmetric positive definite matrices $P,Q$ such that 
\bea A^TP+PA&=& -Q\nonumber\\ PB=C\label{e:kyl}\eea

\subsection*{Bregman Divergence}
Further extensions can be obtained based on the Bregman divergence construction \cite{Bof21}. Let $f(x)$ be a twice differentiable function, $x\in R^n$. Let $D_f(x,y)=f(x)-f(y)-(\nabla f(y),x-y)$, where $\nabla f(x)$ is the gradient of the function $f(x)$. The function $D_f(x,y)$ turns out to be convenient to use for the convergence proofs as a part of extended Lyapunov function. For example, extended versions of the speed-gradient algorithm can be designed  via a Lyapunov function \cite{Bof21}
\begin{equation} \label{F5-10}
V(x,\theta,t)=Q(x,t)+D_f(\theta,\theta_*).
\end{equation}
It leads to the algorithms
\begin{equation} \label{F5-11}
\dot\theta=-\nabla^2D_f(\theta,\theta_*)\nabla_{\theta} w(x,\theta,t),
\end{equation}
generalizing algorithms (\ref{F5-5}).

\subsection*{Averaging}
A standard method that has been studied extensively in the area of nonlinear oscillations has been utilized often in adaptive systems in the context of robustness with respect to disturbances and unmodeled dynamics \cite{And86}. This is briefly summarized below.
Originally suggested in \cite{Kry43} and expanded in \cite{Bog61,Sanders85}, this method is associated with the solutions of a differential equation
\be \dot x = \mu f(x,t,\mu),\qquad x(0)=x_0\label{average1}\ee
where $\mu$ is a positive constant. An approximate solution for \eqref{average1} can be found if $\dot x$ is small and the solution $x(t)$ varies slowly using the process of averaging. The underlying idea here is that as $x(t)$ is varying slowly, the rapidly varying terms in $f$ do not affect the slow variation of $x$ in the long run. We briefly outline the application of this tool to adaptive systems \cite{And86,Kok85}:
The underlying error model for the perturbed system in \eqref{plant-linear-robust}, when only unmodeled dynamics are present and the adaptive law as in \eqref{ad20} and \eqref{ad21} is used can be written as 
\be
\left[\begin{array}{l} \dot e\\\dot{\wt}\end{array}\right]=\left[\begin{array}{cc}A & b\omega^T\\ -\mu\omega c^T & 0\end{array}\right]\label{average2}\ee
where $\{c,A,b\}$ corresponds to the realization of the closed-loop transfer function that arises when there is no parametric uncertainty, and $\omega$ denotes the system variables that are accessible. That is, $\bar W_m(s)=c^T(sI-A)^{-1}b$. The robustness of the adaptive system is assured if the solutions of \eqref{average2} are well behaved. The following theorem outlines this result \cite{Kok85}:
Let $\omega(t)$ be bounded, almost periodic, and persistently exciting. Then 
\begin{enumerate}
\item there exists a $c^*>0$ such that for all $\mu\in (0,c^*]$, the origin of \eqref{average2} is exponentially stable if 
\be Re\left[\lambda_i\left(\int_0^T\omega(t)\bar W_m(s)\omega^T(t)dt \right)\right] > 0 \;\; \forall i=1,\ldots n \label{average3}\ee
\item The condition \eqref{average3} is satisfied if 
\be \sum_{k=-\infty}^\infty Re\left[\bar W_m(i\nu_k)\right]Re\left[\Omega(i\nu_k)\bar\Omega^T(i\nu_k)\right]>0\label{average4}\ee
\end{enumerate}
An expansion of $\omega(t)$ using an inverse Fourier transform expansion such as $\omega(t)=\sum_{k=-\infty}^\infty \Omega(i\nu_k)exp(i\nu_k t)$ is leveraged in this context. Eq. \eqref{average4} implies that the stability property of \eqref{average2} critically depends on the spectrum of the excitation of $\omega$ in relation to the closed-loop transfer function $\bar W_m(s)$. More importantly, the condition in \eqref{average3} can be met by a large class of problems where $\bar W_m(s)$ is not SPR.


\bibliographystyle{elsarticle-harv}
\bibliography{IEEEabrv, references, References_Fradkov, ALF_BIB-2}
\end{document}